\documentclass[11pt]{article}
\usepackage[margin=1.3in]{geometry}
\usepackage{graphicx} % Required for inserting images
\usepackage{authblk}  % For better author formatting
\usepackage{hyperref} % For clickable email links
\usepackage[dvipsnames]{xcolor}
\usepackage{graphicx}      
\usepackage{amsmath,amssymb, amsfonts,amsthm}
\usepackage{subcaption}
\usepackage{graphicx}
\usepackage[dvipsnames]{xcolor}
\usepackage{enumitem}
\usepackage{url}
\usepackage{caption}
\usepackage{amsmath}
\usepackage{array} % Required for tabular
\usepackage{amsfonts}
\usepackage{amssymb}
\usepackage{nicefrac}
\usepackage{graphicx}      % include this line if your document contains figures
\usepackage{tikz}
\usepackage{pgfplots}
\pgfplotsset{compat=1.8}
\usepackage{pgfplotstable}
\usetikzlibrary{3d, calc, decorations.markings}
\usepackage{placeins}    % Ensures figures do not float past barriers
\pgfplotsset{compat=1.17}

\usepackage[style= numeric-comp,hyperref=true, doi=false,url=false,
            isbn=false,
            firstinits=true, sorting = none, 
            block=none, backend=bibtex,maxnames=99]{biblatex}
\renewbibmacro{in:}{} 
\bibliography{references}

\newtheorem{definition}{Definition}
\newtheorem{theorem}{Theorem}
\newtheorem{lemma}{Lemma}
\newtheorem{proposition}[theorem]{Proposition}
\newtheorem{corollary}[theorem]{Corollary}

\newcommand{\rk}{\operatorname{rank}}
\newcommand{\supp}{\operatorname{supp}}

\newcommand{\N}{\mathbb{N}}

\newcommand{\cU}{\mathsf{U}}
\newcommand{\cA}{\mathsf{A}}
\newcommand{\cX}{\mathsf{X}}

\newcommand{\cE}{\mathcal{E}}
\newcommand{\R}{\mathbb{R}}

\newcommand{\cH}{\mathcal{H}}
\newcommand{\bP}{\mathbf{P}}
\newcommand{\bfo}{\mathbf{1}}
\newcommand{\rint}{\operatorname{int}}
\newcommand{\rs}{\mathsf{s}}
\newcommand{\aas}{{\em a.a.s.}}
\newcommand{\crk}{\operatorname{co-rank}}

\newcommand{\xc}[1]{\vspace{.1cm}

\noindent {\em #1} }

\makeatletter
\def\blfootnote{\xdef\@thefnmark{}\@footnotetext}
\makeatother

\makeatletter
\newcommand{\oset}[3][0ex]{%
  \mathrel{\mathop{#3}\limits^{
    \vbox to#1{\kern-1\ex@
    \hbox{$\scriptstyle#2$}\vss}}}}
\makeatother

\title{Hamiltonicity of Step-graphons}

\begin{document}
\author{Xudong Chen$^\dagger$
}
\date{}
\maketitle
\blfootnote{$^\dagger$X. Chen is with the Department of Electrical and Systems Engineering, Washington University, St. Louis. Email: \texttt{cxudong@wustl.edu}.}

\begin{abstract} 
A step-graphon has the strong (resp., weak) $H$-property if a directed, random graph sampled from it has a Hamilton cycle (resp., a node-wise disjoint cycle cover) asymptotically almost surely. The weak/strong $H$-property is essentially a zero-one property. We identify key objects associated with the step-graphon that matter for the zero-one law and provide a complete characterization.         
\end{abstract}

\section{Introduction}\label{sec:introduction}
In this paper, a graphon $W$ is a measurable function $W: [0,1]^2 \to [0,1]$, which is {\it not} necessarily symmetric. 
We treat graphons as stochastic models and investigate their hamiltonicity. Specifically, given a graphon $W$, we obtain a digraph $\vec G_n \sim W$ on $n$ nodes via the following two-step sampling procedure:
\begin{enumerate}
    \item[$\mathbf{S1}$.] Sample $t_1, \ldots, t_n \sim \text{Uni}[0,1]$ independently, where $\text{Uni}[0,1]$ is the uniform distribution over the interval $[0,1]$. We call $t_i$ the \emph{coordinate} of node $v_i$.
    \item[$\mathbf{S2}$.] For each pair of {\em distinct} nodes $v_i$ and $v_j$, place independently a {\em directed} edge from $v_i$ to $v_j$ with probability $W(t_i, t_j)$ and a {\em directed} edge from $v_j$ to $v_i$ with probability $W(t_j, t_i)$. 
\end{enumerate}

A digraph $\vec G$ is said to have a {\it node-wise disjoint cycle cover} (or, simply, cycle cover) if it contains a subgraph $\vec H$, with the same node set as $\vec G$, such that $\vec H$ is a node-wise disjoint union of directed cycles of $\vec G$. If $\vec H$ is a cycle, then $\vec H$ is a {\em Hamilton cycle} of $\vec G$. We evaluate the probability that $\vec G_n \sim W$ has a cycle cover or a Hamilton cycle as $n\to\infty$. A precise problem formulation will be given shortly.   

Our interest in cycle cover is rooted in structural system theory, which deals with the problem of characterizing what type of network topology can sustain a desired system property. To elaborate, consider a network of $n$ mobile agents whose communication topology is described by a digraph $G$, where the nodes represent the agents and the edges indicate the information flows.    
More specifically, a directed edge from $v_j$ to $v_i$ indicates that agent $x_i$ can access the state information of agent $x_j$, so the dynamics of $x_i$ are allowed to depend on the state of $x_j$. Said in other words, if the dynamics of the network system obey the differential equation $\dot x_i = f_i(x(t))$, for all $i = 1,\ldots, n$, with each $f_i$ a differentiable function, then $\partial f_i(x)/\partial x_j \neq 0$ only if there is an edge from $v_j$ to $v_i$ in $G$. We call such dynamics {\it compatible with $G$}. The central question of the structural system theory is then the following: Given the digraph $G$ and given a desired system property (e.g., asymptotic stability with respect to the origin), is there a dynamical system $\dot x = f(x)$ such that it is compatible with $G$ and satisfies the  property? This framework can be extended to controlled system $\dot x = f(x, u)$, taking into account the constraint that for each control input $u_j$, there may be only few agents $x_i$ under its direct influence, i.e., $\partial f_i(x)/\partial u_j \neq 0$.         

It has been shown that existence of a cycle cover is essential for a network topology to sustain (ensemble) controllability~\cite{chen2021sparse} and  stability~\cite{belabbas2013sparse}, two of the fundamental properties of a dynamical control system. 
When a multi-agent system operates in an uncertain and/or adversarial environment,  its network topology becomes a random object. 
We use a graphon $W$ to represent the environment uncertainty and the random digraph $\vec G_n \sim W$ to represent the network topology, so the probability that an ordered pair of agents establishes an oriented communication link depends on their respective positions.  
The knowledge about how likely the network topology of a large-scale multi-agent system can have a cycle cover is critical for a network manager to understand whether the environment is in favor of or against them, to evaluate the risk-to-reward ratio, and to decide whether the system shall be deployed. 

\xc{Problem formulation.}  We start by introducing the class of step-graphons:

\begin{definition}[Step-graphon]\label{def:stepgraphon}
A graphon $W$ is a {\bf step-graphon} if there is a sequence $0 =: \sigma_0 < \sigma_1 < \cdots < \sigma_m := 1$, for some $m \geq 1$, such that $W$ is constant over the rectangle $R_{ij}:=[\sigma_{i-1}, \sigma_{i}) \times [\sigma_{j-1}, \sigma_{j})$ for  $1\leq i,j \leq m$. We call 
$\sigma := (\sigma_0, \sigma_1, \ldots, \sigma_{m})$ a {\bf partition} for $W$.
\end{definition}
 
We illustrate step-graphons in Figure~\ref{sfig1:stepgraphon}. We next have the following definition: 

\begin{definition}[$H$-property] 
    A graphon $W$ has {\bf weak $H$-property} if $\vec G_n\sim W$ has a node-wise disjoint cycle cover asymptotically almost surely (\aas), i.e.,
    \begin{equation}\label{eq:defHproperty}
    \lim_{n\to\infty}\bP(\vec G_n \sim W \mbox{ has a cycle cover}) = 1. 
    \end{equation} 
    The graphon has {\bf strong $H$-property} if $\vec G_n\sim W$ has a Hamilton cycle \aas, i.e.,
    \begin{equation}\label{eq:defstrongHproperty}
    \lim_{n\to\infty}\bP(\vec G_n \sim W \mbox{ has a Hamilton cycle}) = 1. 
    \end{equation}
\end{definition}

We will see soon that the weak/strong $H$-property is essentially a zero-one property, i.e., for almost all step-graphons, the limit on the left hand side of~\eqref{eq:defHproperty} or~\eqref{eq:defstrongHproperty} is either zero or one.   The main theorem of this paper, which we present in Subsection~\ref{ssec:mainresult}, provides a complete characterization of the zero-one law. 

The notion of weak $H$-property was introduced in the earlier work~\cite{belabbas2021h,belabbas2023geometric}, in which we have established the zero-one law for the class of {\it symmetric} step-graphons. 
Specifically, we have identified a few key objects associated with a step-graphon that matter for the zero-one law and formulated a set of conditions that are necessary and sufficient for the step-graphon to have the weak $H$-property.   
These objects will be modified appropriately to account for the asymmetry of the step-graphon, and will be introduced in Subsection~\ref{ssec:keyobjects}. 
The main theorem of the paper (Theorem~\ref{thm:main}) builds upon these new objects and presents conditions that decide the zero-one law for both weak and strong $H$-property. 

In a more recent work~\cite{chen2025h}, we have addressed the zero-one law for the weak $H$-property of general step-graphons. Specifically, we have provided a set of conditions that are necessary and sufficient for a step-graphon to satisfy the weak $H$-property, together with a sketch of proof.  
The present paper expands significantly on the preliminary version~\cite{chen2025h} by generalizing the result to include the zero-one law for the strong $H$-property and providing a complete proof with finer and more thorough arguments.
   
At the end of this section, we gather a few key notations and terminologies used throughout the paper.

\xc{Notation.} 
In this paper, we consider both directed and undirected graphs. We will put an arrow on top of the letter (e.g., $\vec G$) to indicate that the graph it refers to is directed.  
All graphs considered in the paper do not have multiple edges, but can have self-loops. 
For a graph $\vec G$, let $V(\vec G)$ and $E(\vec G)$ be its node set and edge set, respectively. We use $v_iv_j$ to denote a directed edge from $v_i$ to $v_j$, and use $(v_i,v_j)$ to denote an undirected edge between $v_i$ and $v_j$. A digraph $\vec G$ is said to be {\it strongly connected} if for any two distinct nodes $v_i$ and $v_j$, there exist a path from $v_i$ to $v_j$ and a path from $v_j$ to $v_i$.   

Let $\R_{>0}$ (resp., $\R_{\geq 0}$) be the set of positive (resp., nonnegative) real numbers. Let $\N$ (resp., $\N_0$) be the set of positive (resp., nonnegative) integers. 

Let $\bfo$ be the vector of all ones, and $e_i$ be the $i$th column of the identity matrix. Their dimension will be clear in the context. 
The support of a vector $x$, denoted by $\supp(x)$, is the set of indices~$i$ such that $x_i\neq 0$. Similarly, the support of a matrix $A = [a_{ij}]$, denoted by $\supp(A)$, is the set of indices~$ij$ such that $a_{ij} \neq 0$.  
We will relate the support of a vector (resp., square matrix) to the node (resp., edge) set of a digraph. Specifically, for a vector $x\in \R^n$ and for a digraph $\vec G$ on $n$ nodes, we can treat $\supp(x)$ as a subset of $V(\vec G)$ where $v_i \in \supp(x)$ if and only if $x_i \neq 0$. Similarly, we treat $\supp(A)$ as a subset of $E(\vec G)$ where $v_iv_j\in E(\vec G)$ if and only if $a_{ij} \neq 0$.   
For a subgraph $\vec G'$ of $\vec G$, we let $x|_{\vec G'}$ be the sub-vector of $x$ obtained by deleting any entry $x_i$ such that $v_i\notin \vec G'$.

\section{Main Result}\label{sec:mainresult}

In this section, we identify the key objects associated with a step-graphon that matter for the weak/strong $H$-property, and formulate a set of conditions that  decide whether the step-graphon has the weak/strong $H$-property. We then illustrate and numerically validate the result. Toward the end, we provide a sketch of proof, highlighting the ideas  that will be used to establish the result.     

\subsection{Key objects}\label{ssec:keyobjects}

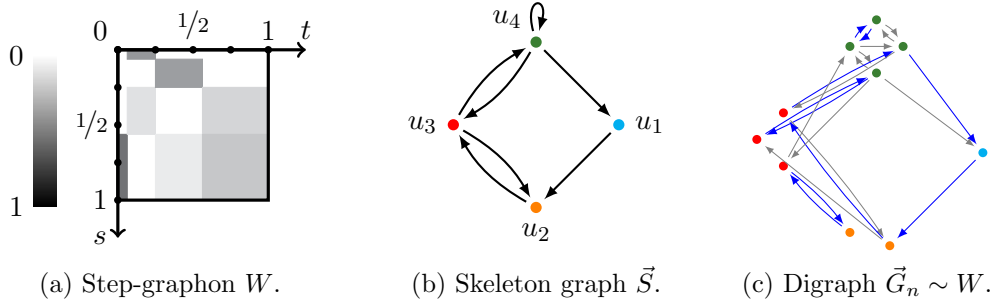
\begin{figure}[t]
    \centering
 \begin{subfigure}{.35\textwidth}
\centering
\begin{tikzpicture}[scale=.5]
    \fill[bottom color=black,top color=white] (0,0) rectangle (0.6,4);
\node [left] at (0,4) {$0$};
\node [left] at (0,0) {$1$};
\node [] at (0,-1) {}; 
\end{tikzpicture}
\begin{tikzpicture}[scale=.5]
\filldraw [fill=Gray!30!black!80, draw=Gray!30!black!80] (0,0) rectangle (1/4,7/4); % 41
\filldraw [fill=Gray!70!black!20, draw=Gray!70!black!20] (1/4,7/4) rectangle (4/4,12/4); % 32
\filldraw [fill=Gray!70!black!70, draw=Gray!70!black!70] (4/4,12/4) rectangle (9/4,15/4); % 23
\filldraw [fill=Gray!10!black!20,draw=Gray!10!black!20] (9/4,7/4) rectangle (16/4,12/4); % 34
\filldraw [fill=Gray!60!black!70, draw=Gray!60!black!70] (1/4,15/4) rectangle (4/4,16/4); % 12
\filldraw [fill=Gray!30!black!10, draw=Gray!30!black!10] (4/4,0) rectangle (9/4,7/4); % 43
\filldraw [fill=Gray!30!black!30, draw=Gray!30!black!30] (9/4,0) rectangle (16/4,7/4); % 44

\draw [draw=black,very thick] (0,0) rectangle (4,4);

\draw[->,very thick,postaction={decorate},
     decoration={
       markings,
       mark = between positions 0 and 0.8 step 0.2 with {\fill[black] circle[radius=0.05cm];}
  }] (0,4) -- (0,-1) node [left] {$s$};
\draw[->,very thick,postaction={decorate},
     decoration={
       markings,
       mark = between positions 0 and 0.8 step 0.2 with {\fill[black] circle[radius=0.05cm];}
  }] (0,4) -- (5,4) node [above] {$t$};
\node [above left] at (0,4) {$0$};
\node [above] at (2,4) {$\nicefrac{1}{2}$};
\node [left] at (0,0) {$1$};
\node [left] at (0,2) {$\nicefrac{1}{2}$};
\node [above] at (4,4) {$1$};
\end{tikzpicture}
\caption{Step-graphon $W$.}\label{sfig1:stepgraphon}
\end{subfigure}
\begin{subfigure}{.3\textwidth}
\centering
 \begin{tikzpicture}[scale=1.1]
 \tikzset{every loop/.style={}}
		\node [circle,fill=OliveGreen,inner sep=1.5pt,label=above left:{$u_4$}] (4) at (0, 0) {};
		\node [circle,fill=red,inner sep=1.5pt,label=left:{$u_3$}] (3) at (-1, -1) {};
		\node [circle,fill=orange,inner sep=1.5pt,label=below:{$u_2$}] (2) at (0, -2) {};
		\node [circle,fill=cyan,inner sep=1.5pt,label=right:{$u_1$}] (1) at (1, -1) {};

  \path[draw,thick,shorten >=2pt,shorten <=2pt, -latex]
		 (1) edge[-latex, black] (2) % 1 2
		 (2) edge[-latex, bend left=15, black] (3) % 2 3
          (4) edge[-latex, black] (1)
          (4) edge[-latex, bend left=15, black] (3)
          (3) edge[-latex, bend left=15, black] (4)
		 (3) edge[-latex, bend left=15, black] (2) % 3 4
		 (4) edge[loop above, -latex, black] (4); % 5 5
\end{tikzpicture}
\caption{Skeleton graph $\vec S$.}\label{sfig2:stepgraphon}
\end{subfigure}
\begin{subfigure}{0.27\textwidth}
\centering
    \begin{tikzpicture}[scale=1,rotate=-45]
    \tikzset{every loop/.style={}}
		\node [circle,fill=red,inner sep=1.2pt] (4) at (0.25, -0.25) {};
		\node [circle,fill=red,inner sep=1.2pt] (5) at (-0.25, 0.25) {}; 
            \node [circle,fill=red,inner sep=1.2pt] (3) at (-0.25, -0.25) {}; 
            %%%% 1
		\node [circle,fill=OliveGreen,inner sep=1.2pt] (6) at (-0.25, 1.5) {};
		\node [circle,fill=OliveGreen,inner sep=1.2pt] (7) at (0.25, 1.5) {};
		\node [circle,fill=OliveGreen,inner sep=1.2pt] (8) at (-0.25, 2) {};
		\node [circle,fill=OliveGreen,inner sep=1.2pt] (9) at (0.25, 2) {}; %%%% 2
		\node [circle,fill=cyan,inner sep=1.2pt] (14) at (2, 1.75) {};
		\node [circle,fill=orange,inner sep=1.2pt] (17) at (1.5, -0.25) {};
		\node [circle,fill=orange,inner sep=1.2pt] (18) at (2, 0) {};
	
  \path[draw, shorten >=2pt,shorten <=2pt, -latex]
         (6) edge[very thin, bend left=15, gray] (7)
         (7) edge[very thin, bend left=15, gray] (6)
		 (6) edge[bend left=15, blue] (8)
         (8) edge[bend left=15, blue] (6)
		 (6) edge[very thin, gray] (9)
		 (9) edge[very thin, bend left=15, gray] (8)
         (8) edge[very thin, bend left=15, gray] (9)
		 (4) edge[bend left=7, blue] (17)
          (17) edge[bend left=7, blue] (4)
		 (3) edge[bend left=4, blue] (7)
          (7) edge[bend left=4, blue] (3)
		 (5) edge[very thin, bend left=4, gray] (18)
          (18) edge[bend left=4, blue] (5)
		 (5) edge[bend left=4, blue] (9)
          (9) edge[very thin, bend left=4, gray] (5)
		 (4) edge[very thin, gray] (6)
		 (7) edge[very thin, gray] (4)
          (18) edge[very thin, gray] (3)
		 (14) edge[blue] (18)
		 (9) edge[blue] (14)
		 (7) edge[very thin, gray] (14)
		 ;
\node at(0,-.51){};
\end{tikzpicture}
\caption{Digraph $\vec G_{n} \sim W$.}\label{sfig2:samplegraph}
\end{subfigure}
\caption{The step-graphon $W$ in~(a) has the partition sequence $\sigma=\frac{1}{16}(0,1,4,9,16)$. The value of $W$ is shade coded, with black being~$1$ and white being $0$.  The digraph $\vec S$ in~(b) is the skeleton graph associated with $W$ with respect to the partition $\sigma$. The digraph $\vec G_{n}$, with $n = 10$, in~(c) is sampled from $W$. It has a cycle cover, highlighted in blue, which comprises a $4$-cycle and three $2$-cycles.}
    \label{fig:stepgraphon}
\end{figure}

We introduce below four objects that are essential to deciding whether a step-graphon has the $H$-property. We start by introducing the definitions of concentration vector and of skeleton graph, which were introduced in~\cite{belabbas2021h,belabbas2023geometric} for symmetric step-graphons and admit natural extensions:  

\begin{definition}[Concentration vector]\label{def:concentrationv} Let $W$ be a step-graphon with partition $(\sigma_0, \ldots, \sigma_m)$. The associated {\bf concentration vector} $x^* = (x^*_1, \ldots, x^*_m)$ has entries defined as 
$x^*_i := \sigma_i - \sigma_{i-1}$, for all $i = 1,\ldots,m$.
\end{definition}

Next, we have 

\begin{definition}[Skeleton graph]\label{def:SkeletonGraph} 
To a step-graphon $W$ with a partition $\sigma = (\sigma_0, \ldots, \sigma_m)$, we assign the digraph $\vec S$ on $m$ nodes $\{u_1, \ldots, u_m\}$, whose edge set $E(\vec S)$ is defined as follows: There is a directed edge from $u_i$ to $u_j$ if and only if $W$ is non-zero over $R_{ij}$. We call $\vec S$ the {\bf skeleton graph} of $W$ for the partition $\sigma$. 
\end{definition}

Note that the concentration vector $x^*$ and the skeleton graph $\vec S$ together determine the support of $W$, and vice versa.  

The next two objects are derived from the skeleton graph $\vec S$. One is the node-cycle incidence matrix $Z$: 

\begin{definition}
[Node-cycle incidence vector/matrix]\label{def:incidencevectorandcone}
Let $\vec C_1,\ldots, \vec C_k$ be the cycles of $\vec S$ (a self-loop is a cycle of length~$1$). 
The {\bf node-cycle incidence vector}  of $\vec C_j$ is   
$z_{j} := \sum_{u_i\in V(\vec C_j)} e_i \in \R^m$. 
The {\bf node-cycle incidence matrix} of $\vec S$ is 
$
Z := 
\begin{bmatrix}
z_1 & \cdots & z_k
\end{bmatrix} \in \R^{m\times k}$.   
\end{definition}

The other object is the convex cone spanned by the column vectors of $Z$:  

\begin{definition}[Node-circulation cone]\label{def:nodeflowcone} 
The {\bf node-circulation cone} $\cX$ of $\vec S$ is the convex cone generated by the node-cycle incidence vectors $z_1,\ldots, z_k$, i.e.,
$$
    \cX := \left \{\sum_{j = 1}^k c_j z_j \mid c_j \geq 0 \mbox{ for all } j = 1,\ldots, k\right \}.
$$
\end{definition}

We name $\cX$ the node-circulation cone for the following reason: A map $f:E(\vec S) \to \R_{\geq 0}$ is said to be a {\em circulation} if it satisfies the following balance condition:
\begin{equation}\label{eq:balancecondition}
\sum_{u_k: u_ku_i\in E(\vec S)} f(u_ku_i) = \sum_{u_j: u_iu_j\in E(\vec S)} f(u_iu_j) =: y_i(f), \quad \mbox{for all } u_i\in V(\vec S). 
\end{equation}
It is not hard to see that $\cX$ is the set of vectors $y(f):= (y_1(f),\ldots, y_m(f))$ for all circulations~$f$.   

It is clear that $\dim \cX = \rk(Z) \leq m$. We will see soon that whether $Z$ has full row rank (i.e., rank $m$) is a deciding factor for weak/strong $H$-property of a step-graphon. 
At the end of the subsection, we present a result that relates the rank of $Z$ to some graphical condition of $\vec S$ (more precisely, the associated bipartite graph $B_{\vec S}$).

\begin{figure}[t]
\centering
    \begin{subfigure}{.4\textwidth}
    \centering
    \begin{tikzpicture}[scale=1]
    \tikzset{every loop/.style={}}
		\node [circle,fill=black,inner sep=1pt,label=above left:{$u_4$}] (4) at (0, 0) {};
		\node [circle,fill=black,inner sep=1pt,label=left:{$u_3$}] (3) at (-1, -1) {};
		\node [circle,fill=black,inner sep=1pt,label=below:{$u_2$}] (2) at (0, -2) {};
		\node [circle,fill=black,inner sep=1pt,label=right:{$u_1$}] (1) at (1, -1) {};

    \path[draw,thick,shorten >=2pt,shorten <=2pt, -latex]
		 (1) edge[-latex, black] (2) % 1 2
		 (2) edge[-latex, bend left=15, black] (3) % 2 3
          (4) edge[-latex, black] (1)
          (4) edge[-latex, bend left=15, black] (3)
          (3) edge[-latex, bend left=15, black] (4)
		 (3) edge[-latex, bend left=15, black] (2) % 3 4
		 (4) edge[loop above, -latex, black] (4); % 5 5
    \end{tikzpicture}
    \caption{Digraph $\vec S$.}\label{sfig:bipar1}
    \end{subfigure}\quad
    \begin{subfigure}{.4\textwidth}
\centering
\begin{tikzpicture}[scale=1]

		\node [circle,fill=black,inner sep=1pt,label=above:{$u'_1$}] (1) at (0, 1.2) {};
		\node [circle,fill=black,inner sep=1pt,label=above:{$u'_2$}] (2) at (1, 1.2) {};
		\node [circle,fill=black,inner sep=1pt,label=above:{$u'_3$}] (3) at (2, 1.2) {};
		\node [circle,fill=black,inner sep=1pt,label=above:{$u'_4$}] (4) at (3, 1.2) {};
		\node [circle,fill=black,inner sep=1pt,label=below:{$u''_1$}] (11) at (0, -0.5) {};
		\node [circle,fill=black,inner sep=1pt,label=below:{$u''_2$}] (12) at (1, -0.5) {};
		\node [circle,fill=black,inner sep=1pt,label=below:{$u''_3$}] (13) at (2, -0.5) {};
		\node [circle,fill=black,inner sep=1pt,label=below:{$u''_4$}] (14) at (3, -0.5) {};
		 
		 \path[draw,thick,shorten >=2pt,shorten <=2pt]
		 (4) edge [black] (11)
		 (1) edge [black] (12)
		 (2) edge [black] (13)
		 (3) edge [black] (12)
          (3) edge [black] (14)
		 (4) edge [black] (13)
		 (4) edge [black] (14);
\end{tikzpicture}
    \caption{Bipartite graph $B_{\vec S}$.}\label{sfig:bipar2}
    \end{subfigure}
    \caption{The bipartite graph $B_{\vec S}$ in~(b) is associated with the digraph $\vec S$ in~(a).}\label{fig:bipartite}
\end{figure}
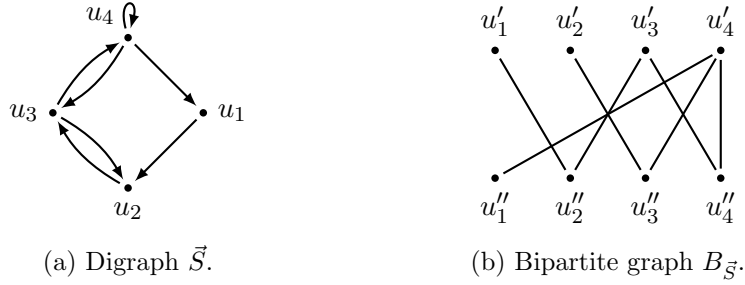

The bipartite graph $B_{\vec S}$ associated with $\vec S$ has $2m$ nodes. The node set $V(B_{\vec S})$ is a disjoint union of two subsets
$V'(B_{\vec S}) := \{u'_1,\ldots, u'_m\}$ and $V''(B_{\vec S}) := \{u''_1,\ldots, u''_m\}$.  
The edge set is given by $E(B_{\vec S}):=\{(u'_i, u''_j) \mid u_iu_j\in E(\vec S)\}$. 
The correspondence between $\vec S$ and $B_{\vec S}$ is illustrated in Figure~\ref{fig:bipartite}.  

Let $\vec S_1, \ldots, \vec S_q$ be the {\it strongly connected components} (SCCs) of $\vec S$. We recall that they satisfy the following three (defining) conditions: (1) Every subgraph $\vec S_p$, for $p = 1,\ldots, q$, is strongly connected; (2) The node sets $V(\vec S_1), \ldots, V(\vec S_q)$ form a partition of $V(\vec S)$; (3) If $\vec S'$ is any other strongly connected subgraph of $\vec S$, then $\vec S'$ is contained in some $\vec S_p$, for $p = 1,\ldots, q$. 

We define the co-rank of $Z$ as $\crk (Z) := m - \rk(Z)$.   
The following result adapted from~\cite[Corollary~5.6]{martens2010separation} provides an explicit formula for the co-rank of~$Z$, which can be evaluated efficiently in $O(|V(\vec S)| + |E(\vec S)|)$ time:

\begin{lemma}\label{lem:connectednessofB}
Let $B_{\vec S_p}$, for $p = 1,\ldots, q$ be the bipartite graph associated with $\vec S_p$. Let $\tau_p$ be the number of connected components of $B_{\vec S_p}$. Then,  
$$
\crk (Z) = \sum_{p = 1}^q (\tau_p - 1). 
$$
In particular, $\crk(Z) = 0$ if and only if every $B_{\vec S_p}$, for $p = 1,\ldots, q$, is connected.  
\end{lemma}

To illustrate, consider the skeleton graph $\vec S$ in Figure~\ref{sfig:bipar1}, with the associated bipartite graph plotted in Figure~\ref{sfig:bipar2}. 
Since $\vec S$ is strongly connected and $B_{\vec S}$ is connected, by Lemma~\ref{lem:connectednessofB} we have that $\crk(Z) = 0$. 
Indeed, $\vec S$ has $4$ cycles: 
\begin{equation}\label{eq:cyclesofexampleS}
\vec C_1 = u_4u_4, \quad \vec C_2 = u_3u_4u_3, \quad \vec C_3 = u_1 u_2 u_3 u_4 u_1, \quad \mbox{and} \quad \vec C_4 = u_2u_3u_2.
\end{equation} 
Thus, one can obtain the node-cycle incidence matrix $Z$ explicitly as 
\begin{equation}\label{eq:defmatrixZ}
Z = 
\begin{bmatrix}
0 & 0 & 1 & 0 \\
0 & 0 & 1 & 1 \\
0 & 1 & 1 & 1 \\
1 & 1 & 1 & 0
\end{bmatrix}, \quad \mbox{so } \crk(Z) = 0.
\end{equation}

\subsection{Zero-one law for $H$-property}\label{ssec:mainresult} 
Let $W$ be a step-graphon and $\sigma$ be a partition for $W$. Let $x^*$, $S$, $Z$, and $\cX$ be the associated concentration vector, skeleton graph, node-cycle incidence matrix, and node-circulation cone, respectively.  
We introduce the following conditions for the quadruple $(x^*,S,Z,\cX)$: 

\begin{description}
    \item[\it Condition $A$:] $\crk(Z) = 0$. 
    \item[\it Condition $B$:] $x^*\in \rint \cX$, where $\rint \cX$ stands for the {\it relative interior} of $\cX$.
    \item[\it Condition $B'$:] $x^*\in \cX$. 
    \item[\it Condition $C$:] $\vec S$ is strongly connected.
\end{description}

These four conditions, though stated with respect to a specific~$\sigma$, are in fact independent of the choice of a partition. Precisely, we have 

\begin{proposition}\label{prop:equivalenceofpartition}
Let $W$ be a step-graphon. Let $\sigma$ and $\sigma'$ be two partitions for $W$, and $(x^*,S,Z,\cX)$ and $(x'^*,S',Z',\cX')$ be the associated quadruples. 
Then, the following items hold: 
\begin{enumerate}
\item Suppose that both $\vec S$ and $\vec S'$ have at least two nodes; then, 
$\vec S$ is strongly connected if and only if $\vec S'$ is.
\item $\crk(Z) = 0$ if and only if $\crk(Z') = 0$.  
\item $x^*\in \cX$ if and only if $x'^*\in \cX'$ ($x^*\in \rint \cX$ if and only if $x'^*\in \rint \cX'$). 
\end{enumerate}
\end{proposition}

The hypothesis of item 1 (i.e., both $\vec S$ and $\vec S'$ have at least two nodes) is meant to rule out the case where $W$ is the zero function. 
To wit, if $W = 0$ and if the partition $\sigma$ is chosen such that $\sigma = (0,1)$, then the associated skeleton graph $\vec S$ comprises a single node $u$ without self-loop. By default, $\vec S$ is strongly connected. However, any other partition $\sigma'$ for $W$ gives rise to a skeleton graph $\vec S'$ such that $\vec S$ has multiple nodes but without any edge.

We provide a proof of the above Proposition~\ref{prop:equivalenceofpartition} in Appendix~\ref{sec:partitionsame}.  
With the result, we can now have the following definition:

\begin{definition}
    We say that a step-graphon $W$ satisfies Condition~$\star$, for some $\star \in\{  A, B, B', C\}$, if there is a partition $\sigma$ for $W$, with $|\sigma| \geq 2$, such that the associated quadruple $(x^*, \vec S, Z, \cX)$ satisfies Condition~$\star$.     
\end{definition}

We can now state the main result of the paper:

\begin{theorem}[Main Theorem]\label{thm:main}
    Let $W$ be a step-graphon. The following hold:
    \begin{enumerate}
    \item If $W$ does not satisfy Condition~$A$ or~$B'$, then 
    \begin{equation}\label{eq:nonevent}
    \lim_{n\to\infty} \bP(\vec G_n \sim W \mbox{ has a cycle cover}) = 0.
    \end{equation} 
    \item If $W$ satisfies Conditions~$A$ and~$B$, but not~$C$, then
    \begin{equation}\label{eq:posevent}
    \lim_{n\to\infty} \bP(\vec G_n \sim W \mbox{ has a cycle cover}) = 1,
    \end{equation} 
    and 
    \begin{equation}\label{eq:posnonevent}
    \lim_{n\to\infty} \bP(\vec G_n \sim W \mbox{ has a Hamilton cycle}) = 0.
    \end{equation} 
    \item If $W$ satisfies Conditions~$A$, $B$, and $C$, then
    \begin{equation}\label{eq:posevent1}
    \lim_{n\to\infty} \bP(\vec G_n \sim W \mbox{ has a Hamilton cycle}) = 1.
    \end{equation} 
    \end{enumerate}
\end{theorem}

As mentioned earlier, the above theorem extends the results of~\cite{belabbas2021h,belabbas2023geometric}. 
We substantiate  our claim in Appendix~\ref{sec:symmetricgraphon}.     

\subsection{Illustration and numerical validation}

\begin{figure}[h]
\centering
\begin{minipage}{.3\textwidth}
\centering
\subfloat[$W_a$]{
\centering
\begin{tikzpicture}[scale=.38]
\filldraw [fill=Gray!50!black!30, draw=Gray!50!black!30] (0,0) rectangle (1/4,7/4); % 41
\filldraw [fill=Gray!50!black!30, draw=Gray!50!black!30] (1/4,7/4) rectangle (4/4,12/4); % 32
\filldraw [fill=Gray!50!black!30, draw=Gray!50!black!30] (4/4,12/4) rectangle (9/4,15/4); % 23
\filldraw [fill=Gray!50!black!30, draw=Gray!50!black!30] (9/4,7/4) rectangle (16/4,12/4); % 34
\filldraw [fill=Gray!50!black!30, draw=Gray!50!black!30] (1/4,15/4) rectangle (4/4,16/4); % 12
\filldraw [fill=Gray!50!black!30, draw=Gray!50!black!30] (4/4,0) rectangle (9/4,7/4); % 43
\filldraw [fill=Gray!50!black!30, draw=Gray!50!black!30] (9/4,0) rectangle (16/4,7/4); % 44
\draw [draw=black,very thick] (0,0) rectangle (4,4);
\end{tikzpicture}
}\,
\subfloat[$W_b$]{
\centering
\begin{tikzpicture}[scale=.38]
\filldraw [fill=black, draw=black] (0,0) rectangle (4/8,8/8); % 41
\filldraw [fill=black, draw=black] (4/8,8/8) rectangle (12/8,20/8); % 32
\filldraw [fill=black, draw=black] (12/8,20/8) rectangle (24/8,28/8); % 23
\filldraw [fill=black, draw=black] (24/8,8/8) rectangle (16/4,20/8); % 34
\filldraw [fill=black, draw=black] (4/8,28/8) rectangle (12/8,16/4); % 12
\filldraw [fill=black, draw=black] (12/8,0) rectangle (24/8,8/8); % 43
\filldraw [fill=black, draw=black] (24/8,0) rectangle (32/8,8/8); % 44

%% x^* = (1,2,3,2)/8

\draw [draw=black,very thick] (0,0) rectangle (4,4);
\end{tikzpicture}
}\vspace{.2cm}

\subfloat[$W_c$]{
\centering
\begin{tikzpicture}[scale=.38]
\filldraw [fill=black, draw=black] (0,0) rectangle (1,0.8); % 41
\filldraw [fill=black, draw=black] (1,0.8) rectangle (2,2); % 32
\filldraw [fill=black, draw=black] (2,2) rectangle (3.2,3); % 23
\filldraw [fill=black, draw=black] (3.2,0.8) rectangle (4,2); % 34
\filldraw [fill=black, draw=black] (1,3) rectangle (2,4); % 12
\filldraw [fill=black, draw=black] (2,0) rectangle (3.2,0.8); % 43
\filldraw [fill=black, draw=black] (3.2,0) rectangle (4,0.8); % 44

%x^* = (1,1,1.2,0.8)/4

\draw [draw=black,very thick] (0,0) rectangle (4,4);
\end{tikzpicture}
}\,
\subfloat[$W_d$]{
\centering
\begin{tikzpicture}[scale=.38]
\filldraw [fill=black, draw=black] (0,0) rectangle (4/8,8/8); % 41
\filldraw [fill=black, draw=black] (4/8,8/8) rectangle (12/8,20/8); % 32
\filldraw [fill=black, draw=black] (12/8,20/8) rectangle (24/8,28/8); % 23
\filldraw [fill=black, draw=black] (24/8,8/8) rectangle (16/4,20/8); % 34
\filldraw [fill=black, draw=black] (4/8,28/8) rectangle (12/8,16/4); % 12
\filldraw [fill=black, draw=black] (12/8,0) rectangle (24/8,8/8); % 43

%% x^* = (1,2,3,2)/8

\draw [draw=black,very thick] (0,0) rectangle (4,4);
\end{tikzpicture}
}

\subfloat[$\vec S$]{
\begin{tikzpicture}[scale=1]
   \tikzset{every loop/.style={min distance = 4mm}}
		\node [circle,fill=black,inner sep=1.5pt,label=above:{$u_4$}] (4) at (0, 1) {};
		\node [circle,fill=black,inner sep=1.5pt,label=below:{$u_3$}] (3) at (0, 0) {};
		\node [circle,fill=black,inner sep=1.5pt,label=below:{$u_2$}] (2) at (1, 0) {};
		\node [circle,fill=black,inner sep=1.5pt,label=above:{$u_1$}] (1) at (1, 1) {};

  \path[draw,thick,shorten >=2pt,shorten <=2pt, -latex]
		 (1) edge[-latex, black] (2) % 1 2
		 (2) edge[-latex, bend left=15, black] (3) % 2 3
          (4) edge[-latex, black] (1)
          (4) edge[-latex, bend left=15, black] (3)
          (3) edge[-latex, bend left=15, black] (4)
		 (3) edge[-latex, bend left=15, black] (2) % 3 4
		 (4) edge[loop left, -latex, black] (4); % 5 5
\end{tikzpicture}
}\,
\subfloat[$\vec S'$]{
\begin{tikzpicture}[scale=1]
   \tikzset{every loop/.style={min distance = 6mm}}
		\node [circle,fill=black,inner sep=1.5pt,label=above:{$u_4$}] (4) at (0, 1) {};
		\node [circle,fill=black,inner sep=1.5pt,label=below:{$u_3$}] (3) at (0, 0) {};
		\node [circle,fill=black,inner sep=1.5pt,label=below:{$u_2$}] (2) at (1, 0) {};
		\node [circle,fill=black,inner sep=1.5pt,label=above:{$u_1$}] (1) at (1, 1) {};

  \path[draw,thick,shorten >=2pt,shorten <=2pt, -latex]
		 (1) edge[-latex, black] (2) % 1 2
		 (2) edge[-latex, bend left=15, black] (3) % 2 3
          (4) edge[-latex, black] (1)
          (4) edge[-latex, bend left=15, black] (3)
          (3) edge[-latex, bend left=15, black] (4)
		 (3) edge[-latex, bend left=15, black] (2) % 3 4
          ;
\end{tikzpicture}
}
\end{minipage}
\begin{minipage}{.69\textwidth}
\centering
\begin{tikzpicture}[scale=1]
\begin{semilogxaxis}[
	xlabel={\small $n$},
	ylabel={\small $p(n)$},
	ymax=1,
	ymin=0
]
% Case (a) x^* = (1,3,5,7)/16, p = 0.2

\addplot[mark=star,mark size=2pt,color=black,thick] coordinates {
	(10, 0.0006) (50, 0.313)  (100,0.8367)
	(500,0.99995)  (1000,1.0) (2000,1.0) (5000, 1.0)
};
%  Case (b) x^* = (1,2,3,2)/8, p=1 with self-loop
\addplot[mark=square*,mark size=1.5pt,color=red,thick] coordinates {
	(10, 0.353) (50, 0.4442)  (100,0.4867)
	(500,0.51795)  (1000,0.51095) (2000,0.50995) (5000,0.50235) 
};

%   Case (c) x^* = (10, 10, 12, 8)/40 p = 1
\addplot[mark=diamond*, mark options={fill},color=brown, mark size=2.5pt,thick] coordinates {
	(10, 0.15385) (50, 0.08365)  (100,0.0536)
	(500,0.00425)  (1000,0.0003) (2000,0.0) (5000,0.0) 
};
%  Case (d)  x^* = (1,2,3,2)/8, p = 1
% WITHOUT self-loop u4
\addplot[mark=star,mark size=2pt,color=blue,thick] coordinates {
	(10, 0.1918) (50, 0.10795)  (100, 0.0779)
	(500,0.03615)  (1000,0.0255) (2000,0.01705) (5000,0.0106)
};
\legend{$(a)$, $(b)$, $(c)$, $(d)$}
\end{semilogxaxis}
\end{tikzpicture}
\end{minipage}
\caption{{\it Left:} Four step-graphons and the associated skeleton graphs, where $\vec S$ in (e) corresponds to $W_a$, $W_b$, $W_c$, and $\vec S'$ in (f) corresponds to $W_d$.  {\it Right:} The empirical probability $p(n)$ that $\vec G_n\sim W_\star$ has a cycle cover, with $20,000$ samples for each $n = 10,50,100,500,1000,2000,5000$.}\label{fig:case4}
\end{figure}
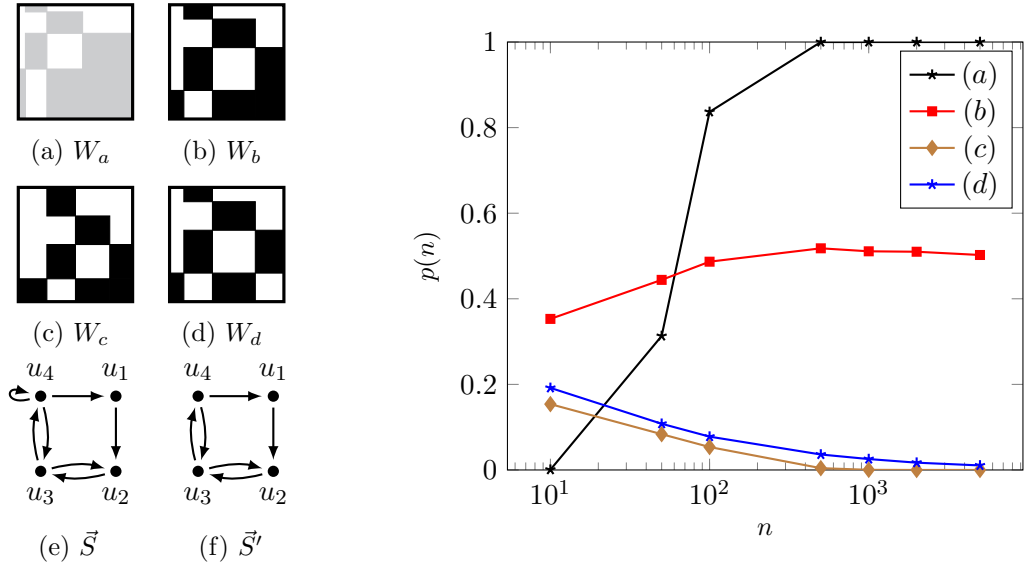

To illustrate Theorem~\ref{thm:main} (more precisely, the zero-law for the weak $H$-property), we consider the four step-graphons in Figure~\ref{fig:case4}. Over their respective support, $W_a$ takes value $0.2$ while $W_b$, $W_c$, and $W_d$ take value $1$.   
We let the partitions for the four step-graphons be 
$$
\sigma_a  = \frac{1}{16}(0,1,4,9,16), \quad  
\sigma_c  = \frac{1}{20}(0,5,10,16,20), \quad \sigma_b =  \sigma_d = \frac{1}{8}(0,1,3,6,8).
$$
The step-graphons in (a), (b), (c) share the same skeleton graph $\vec S$ as shown in (e), which is the same as the one in Figure~\ref{sfig2:stepgraphon}. 
The skeleton graph $\vec S'$ associated with the step-graphon in (d) is shown in (f), which can be obtained from $\vec S$ by removing the self-loop $u_4u_4$. 

The skeleton graph $\vec S$ has $4$ cycles $\vec C_1,\ldots, \vec C_4$ as shown in~\eqref{eq:cyclesofexampleS}. The node-cycle incidence matrix $Z$ has full row rank as shown in~\eqref{eq:defmatrixZ}.  
The digraph $\vec S'$, being a subgraph of $\vec S$, has only three cycles $\vec C_2, \vec C_3, \vec C_4$. Its node-cycle incidence matrix $Z'$ is given by
\begin{equation}\label{eq:defmatrixZ'}
Z' := 
\begin{bmatrix}
0 & 1 & 0 \\
0 & 1 & 1 \\
1 & 1 & 1 \\
1 & 1 & 0
\end{bmatrix}, \quad \mbox{so }\crk(Z') = 1.
\end{equation}

We state without a proof that  any three column vectors of $Z$ form a facet-defining hyperplane of the cone $\cX$. For each $i = 1,\ldots,4$, we let $L_i$ be the subspace spanned by the $z_j$'s, for $j\neq i$. Let $g_i\in \R^4$ be the normal vector perpendicular to $L_i$ of unit length such that $g_i^\top z_i > 0$. Then, it is not hard to obtain that 
$$
g_1 = \frac{1}{2}(-1,1,-1,1), \,\, g_2 = \frac{1}{\sqrt{2}}(0,-1,1,0), \,\, g_3 = (1,0,0,0), \,\, g_4 = \frac{1}{\sqrt{2}}(-1,1,0,0).
$$
Using the half-space representation, we write
$$
\cX = \{y\in \R^4 \mid g_i^\top y \geq 0 \mbox{ for all } i = 1,\ldots, 4\}. 
$$

We numerically validate the necessity and sufficiency of Conditions $A$, $B$ (or $B'$) for the step-graphon $W_\star$, for $\star = a,b,c,d$, to have the weak $H$-property. For each case and for each $n \in \{10, 50, 100, 500, 1000, 2000, 5000\}$, we sampled $20,000$ random graphs $\vec G_n \sim W_\star$ and plot the empirical probability $p(n)$ that $\vec G_n$ has a cycle cover, i.e., 
$$p(n):= \frac{\mbox{number of } \vec G_n\sim W \mbox{ has a cycle cover}}{20,000}.$$

\xc{Case (a).} The concentration vector is $x^*_a = \frac{1}{16}(1,3,5,7)$. 
It belongs to the relative interior of $\cX$ as we can express $x^*_a$ as a positive combination of the $z_j$'s (the column vectors of $Z$ given in~\eqref{eq:defmatrixZ}):  
$$
x^*_a =  \frac{1}{4} z_1 + \frac{1}{8} z_2 + \frac{1}{16} z_3 + \frac{1}{8} z_4. 
$$
Also, the matrix $Z$ has full row rank. Thus, $W_a$ satisfies Conditions $A$ and $B$.  We see from the simulation that the empirical probability $p(n)$ converges to $1$. 

\xc{Case (b).} The concentration vector is $x^*_b = \frac{1}{8}(1,2,3,2)$. It belongs to the {\it boundary} of $\cX$, i.e., $x^*_b \in \cX - \rint \cX$. To wit, we note that $x^*_b$ is in the facet-defining hyperplane $L_1$ spanned by $z_2$, $z_3$, and $z_4$:
$$
x^*_b = \frac{1}{8} (z_2 +  z_3 + z_4).
$$
Thus, $W_b$ satisfies Conditions $A$ and $B'$, but not $B$. We see from the simulation that the empirical probability $p(n)$ converges to neither $1$ nor $0$. In fact, using the same arguments as in~\cite{GAO20257}, we can show that the probability converges to $0.5$.  

\xc{Case (c).} The concentration vector is $x^*_c = \frac{1}{20}(5,5,6,4)$.  
Since $g_1^\top x^*_c < 0$, $x^*_c$ does not belong to $\cX$.  
Thus, $W_c$ satisfies Condition $A$ but not $B'$. 
We see from the simulation that the empirical probability $p(n)$ converges to $0$.

\xc{Case (d).} The concentration vector is $x^*_d = \frac{1}{8}(1,2,3,2)$, same as the one in Case (b). As argued above, we can write $x^*_d$ as a positive combination of $z_2$, $z_3$, and $z_4$, which are the three column vectors in $Z'$. Thus, $x^*_d\in \rint \cX'$. Also, as shown in~\eqref{eq:defmatrixZ'}, $Z'$ does not have full row rank. Thus, $W_d$ satisfies Condition $B$ but not $A$.   
We see from the simulation that the empirical probability $p(n)$ converges to $0$. 

\subsection{Sketch of proof of Theorem~\ref{thm:main}}
We start by introducing two new objects, which will be of great use in the proof. 

\begin{definition}[$\vec S$-partite graph]\label{def:spartite}
    Let $\vec S$ be an arbitrary digraph on $m$ nodes, possibly with self-loops. 
    A directed graph $\vec G$ is {\bf $\vec S$-partite} if there exists a graph homomorphism $\pi:\vec G \to \vec S$.  
    Further, $\vec G$ is a {\bf complete $\vec S$-partite graph} if $$v_iv_j \in E(\vec G)  \quad \Longleftrightarrow \quad \pi(v_i)\pi(v_j) \in E(\vec S).$$ 
\end{definition}

For an $\vec S$-partite graph $\vec G$, we let  
$$y(\vec G) :=  (y_1, \ldots, y_m) \quad \mbox{with } y_i:= |\pi^{-1}(u_i)|, \quad \mbox{for all } i = 1,\ldots, m.$$ 
Further, for a vector $y \in \N^m_0$, let $\vec K_{y}$ be the complete $\vec S$-partite graph, with $y(\vec K_{y}) = y$. In case we need to emphasize the dependence of $\vec K_y$ on $\vec S$, we write $\vec K_y(\vec S)$. 

The relevance of $\vec S$-partite graphs is apparent. Any random graph $\vec G_n$ sampled from $W$ is $\vec S$-partite, where the homomorphism $\pi: \vec G_n \to \vec S$ is naturally the one that sends each node $v_j\in \vec G_n$, with coordinate $t_j\in [\sigma_{i-1}, \sigma_i)$, to $u_i$. 
It is also clear from the sampling procedure (more specifically, step $S1$) that $y(G_n)$ is a multinomial random variable with $n$ trials, $m$ events, and $x^*_i$'s the event probabilities.  
Let 
\begin{equation}\label{eq:defxGn}
x(\vec G_n) := \frac{1}{n}y(\vec G_n). 
\end{equation}
We call $x(\vec G_n)$ the {\it empirical concentration vector}. It follows directly from the (strong) law of large numbers that  
\begin{equation}\label{eq:empiricalx}
x(\vec G_n) \to x^* \mbox{ \aas}.
\end{equation}

Next, we have the following definition: 

\begin{definition}[Edge-circulation cone]\label{def:CaS}
To the skeleton graph $\vec S$ on $m$ nodes,  we assign the set $\cA$ of $m\times m$ nonnegative matrices $A$ such that $\supp(A)\subseteq E(\vec S)$ and 
\begin{equation}\label{eq:balanceequationforA}
A^\top \bfo = A \bfo.
\end{equation}  
We call $\cA$ the {\bf edge-circulation cone}. 
\end{definition}

It is clear from the definition that $\cA$ is a convex cone. Its relation to circulations is as follows: 
Let $f: E(\vec S)\to \R_{\geq 0}$ be a circulation, and $A(f) = [a_{ij}(f)]$ be  such that 
$$
a_{ij}(f) :=
\begin{cases}
f(u_iu_j) & \mbox{if } u_i u_j \in E(\vec S), \\
0 & \mbox{otherwise}. 
\end{cases}
$$
The balance condition~\eqref{eq:balancecondition} of~$f$ guarantees that~\eqref{eq:balanceequationforA} is satisfied, so $A(f)\in \cA$. It turns out that $\cA$ is the set of $A(f)$ for all circulations $f$ on $\vec S$. 
For each cycle $\vec C_j = u_{j_1}u_{j_2}\cdots u_{j_d}u_{j_1}$ of $\vec S$, let $A_j:= \sum_{i = 1}^{d} e_{j_i}e_{j_{i + 1}}^\top$ be the associated adjacency matrix, where $j_{d + 1}$ is identified with $j_1$. It is clear that $A_j \in \cA$ and that $\cA$ is generated by the $A_j$'s.    
Also, note that 
$z_j = A_j \bfo$, so the two cones $\cX$ and $\cA$ relate to each other in the following way: 
\begin{equation}\label{eq:relationya}
\cX = \cA \bfo = \{A\bfo \mid A \in \cA\}. 
\end{equation}

Our use of the edge-circulation cone is through the cycle covers of $\vec S$-partite graphs. Specifically, if $\vec G$ is an $\vec S$-partite graph and if $\vec G$ has a cycle cover $\vec H$, then $\vec H$ induces an integer valued circulation $f_{\vec H}$ on $\vec S$ in the way such that $f_{\vec H}(u_iu_j)$ records the total number of edges used in $\vec H$ from the nodes in $\pi^{-1}(u_i)$ to the nodes in $\pi^{-1}(u_j)$. Consider, for example, the digraph $\vec G_{n}$ in Figure~\ref{sfig2:samplegraph}, with $n = 10$. It has a cycle cover $\vec H$ highlighted in blue. Then, the corresponding $A$-matrix is given by 
$$A(f_{\vec H}) = 
\begin{bmatrix}
0 & 1 & 0 & 0 \\
0 & 0 & 2 & 0 \\
0 & 1 & 0 & 2 \\
1 & 0 & 1 & 2 
\end{bmatrix}.
$$  
It follows from the construction (see Lemma~\ref{lem:necessaryyinY} for a proof) that $y(\vec G) = A(f_{\vec H}) \bfo \in \cX$.

With the $\vec S$-partite graphs and the edge-circulation cone introduced above, we now sketch the proof of Theorem~\ref{thm:main}:   

\subsubsection{On necessity of Conditions $A$, $B'$, and $C$} 
As argued above, if $\vec G_n\sim W$ has a cycle cover, then $y(\vec G_n)\in \cX$. Thus, to establish item~1 of Theorem~\ref{thm:main}, 
it suffices to show that 
$$\neg A \, (\crk(Z) > 0) \mbox{ or } \neg B' \, (x^*\notin \cX) \quad \Longrightarrow \quad y(\vec G_n)\notin \cX \mbox{ \aas}.$$  
The proof that $$\neg B' \quad \Longrightarrow \quad y(\vec G_n) \notin \cX \mbox{ \aas}$$ 
is straightforward, following directly from~\eqref{eq:empiricalx}. The proof that $$\neg A \quad  \Longrightarrow \quad y(\vec G_n) \notin \cX \mbox{ \aas}$$ uses the following arguments: Let $\Delta^{m-1}$ be the standard simplex  in $\R^m$  
and $\overline \cX:= \cX \cap \Delta^{m-1}$. It is not hard to see that $y(G_n) \in \cX$ if and only if $x(G_n)\in \overline \cX$, where we recall that $x(G_n)$ is the empirical concentration vector~\eqref{eq:defxGn}. 
Note that if $\crk(Z) \geq 1$, then $\dim \overline \cX < \dim \Delta^{m-1} = m - 1$. Appealing to the central limit theorem, we have that the random variable $\omega(G_n) := \sqrt{n}(x(G_n) - x^*) + x^*$ converges in distribution to the Gaussian random variable $\omega^*$ whose support is known to be the entire hyperplane that contains $\Delta^{m-1}$. As a consequence, it holds that $x(G_n)\notin \overline\cX$ \aas.

The necessity of Condition $C$ ($\vec S$ is strongly connected) for the strong $H$-property follows from the fact if $\vec H$ is a Hamilton cycle of $\vec G_n$, then $\pi(\vec H)$ is a closed walk of $\vec S$. It is an immediate consequence of~\eqref{eq:empiricalx} that $\pi(\vec H)$ visits every node of $\vec S$~\aas, and hence, $\vec S$ must be strongly connected. 

A complete proof of the necessity part will be presented in Section~\ref{sec:necessityofAB'}.  

\subsubsection{On sufficiency of Conditions $A$, $B$, and $C$} 
We introduce a subset $\cX_0$ of $\cX$, which comprises all integer-valued $y\in \cX$ such that $\|y\|_1$ is sufficiently large and $y/\|y\|_1$ is sufficiently close to $x^*$. A precise definition of $\cX_0$ will be given 
at the beginning of Section~\ref{sec:constructHam}. 
The two conditions $A$ and $B$, together with~\eqref{eq:empiricalx}, guarantee that $y(G_n)\in \cX_0$ \aas. 
The major task is then to show that
\begin{multline*}\label{eq:task1}
x(G_n) \in \cX_0 \mbox{ (and } \vec S \mbox{ is strongly connected)}   \quad \Longrightarrow \quad \\ G_n \mbox{ has a cycle cover (Hamilton cycle) \aas}. 
\end{multline*}
To accomplish the task, we take a two-step approach: 

\xc{Step~1:} We show that if $y\in \cX_0$ (and if $\vec S$ is strongly connected), then the complete $\vec S$-partite graph $\vec K_y$ has a cycle cover (Hamilton cycle). 
The proof builds upon the following facts:
\begin{enumerate}
\item[\it 1.1.] The first fact is an implication of the integrality theorem for the maximum flow problem, which says that if $y \in \cX$ is {\it integer valued}, then there exists an {\it integer-valued} $A\in \cA$ such that $A\bfo = y$. 
\item[\it 1.2.] We then express the matrix $\cA$, obtained from above, as an integer combination of the adjacency matrices $A_j$ associated with the cycles $\vec C_j$  of $\vec S$, i.e., we write $A = \sum_{j = 1}^k c_j A_j$ for $c_j \in \N_0$ (and $c_j\in \N$ in the case $y\in \cX_0$).  
We show that $\vec K_y$ has a cycle cover, which contains $c_j$ cycles that are isomorphic to $\vec C_j$ under the map $\pi: \vec K_y \to \vec S$.   
\item[\it 1.3.] If, further, $\vec S$ is strongly connected, then the cycles  of the  cycle cover exhibited above can be used to form a Hamilton cycle. The proof will be carried out by induction on the number of cycles in $\vec S$ and relies on the use of the (directed) ear decomposition of $\vec S$.
\end{enumerate}
Complete arguments for this step will be presented in Section~\ref{sec:constructHam}. 

\xc{Step~2:} Let $\vec H_y$ be the cycle cover (Hamilton cycle) of $\vec K_y$ obtained in Step~1. We show that $\vec G_n$ contains $\vec H_{y(\vec G_n)}$ as a subgraph~\aas. Precisely, let $\psi_y: \vec H_y \to \vec K_y$ be the embedding. Composing $\psi_y$ with $\pi$, we obtain the map $\pi\cdot \psi_y: \vec H_y \to \vec S$. 
We show that \aas~there exists an embedding $\phi: \vec H_{y(\vec G_n)} \to \vec G_n$ such that $\phi$ is compatible with $\psi_{y(\vec G_n)}$, i.e., $\pi\cdot \phi = \pi\cdot \psi_{y(\vec G_n)}$. The proof relies on the use of the {\it Blow-up Lemma}~\cite{komlos1997blow}. Roughly speaking, the lemma states that if an {\it undirected} graph $H$ has its degree bounded above by a constant and if it can be embedded into a complete $S$-partite graph $K_y$, where $S$ is an {\it undirected} graph {\it without} self-loop, then the graph $H$ can be embedded into any $S$-partite graph $G$, with $y(G) = y$, as long as $G$ satisfies some regularity condition.  
To enable its use, we take the following steps:
\begin{enumerate}
\item[\it 2.1] In Section~\ref{sec:loopfree}, we show that if the step-graphon $W$ has a nonzero diagonal block (i.e., $\vec S$ has a self-loop) and satisfies Condition~$\star$, for $\star = A, B, C$, then there is a step-graphon $W'$ such that $W'\leq W$ (i.e., $W'(s,t)\leq W(s,t)$ for all $(s,t)\in [0,1]^2$), $W'$ satisfies Condition~$\star$ and, moreover, $W'$ is ``loop free'', i.e., the associated skeleton graph does not have any self-loop.  
This fact, combined with the monotonicity of the (strong) $H$-property, allow us to consider only the class of loop-free step-graphons for establishing the sufficiency of Conditions $A$, $B$, and $C$.

\item[\it 2.2] In Section~\ref{sec:sufficiency}. we introduce 
an auxiliary symmetric step-graphon $W^{\rs}$, which is derived from~$W$, together with an auxiliary sampling procedure that allows us to draw undirected random graphs $G_n$ from $W^{\rs}$. The graphon $W^{\rs}$ and the sampling procedure are defined in a way such that the probability that $\vec H_{y(\vec G_n)}$ is embeddable into $\vec G_n$ is bounded above by the probability that $H_{y(G_n)}$ is embeddable into $G_n$, where $H_{y(G_n)}$ is the undirected counterpart of $\vec H_{y(\vec G_n)}$. 
\end{enumerate}
We then complete the proof by showing that~\aas~the random graph $G_n$ satisfies the aforementioned regularity condition. Thanks to the Blow-up lemma, $H_{y(G_n)}$ can be embedded into $G_n$~\aas.

\section{Proof of the Necessity of Conditions $A$, $B'$, and $C$}\label{sec:necessityofAB'}
In this section, we establish ({\it i}\,) the necessity of Conditions $A$ (i.e., $\crk Z = 0$) and $B'$ (i.e., $x^*\in \cX$) for weak $H$-property, and ({\it ii}\,) the necessity of Condition $C$ (i.e., $\vec S$ is strongly connected) for strong $H$-property. 

We start with the following lemma, which establishes a necessary condition for an $\vec S$-partite graph to have a cycle cover: 

\begin{lemma}\label{lem:necessaryyinY}
Let $\vec G$ be an $\vec S$-partite graph. If $\vec G$ has a cycle cover, then $y(\vec G)\in \cX$. 
\end{lemma}

\begin{proof}
Let $\pi: \vec G\to \vec S$ be the graph homomorphism, and $\vec H$ be a cycle cover of $\vec G$. 
Let $n_{ij}$, for $1\leq i, j \leq m$, be the number of directed edges of $\vec H$ from nodes in $\pi^{-1}(u_i)$ to nodes in $\pi^{-1}(u_j)$. 
It is clear that for all $u_i \in V(\vec S)$, 
\begin{equation}\label{eq:solH}
|\pi^{-1}(u_i)|=\sum_{j = 1}^m n_{ij} = \sum_{j = 1}^m n_{ji}.
\end{equation}
Now, consider the matrix 
$A:= \left [ n_{ij}\right ]_{1\leq i, j\leq m}$. 
It is clear that $\supp(A)\subseteq E(\vec S)$. 
Also, by~\eqref{eq:solH}, we have that 
$$ A^\top \bfo = A\bfo = y(\vec G),$$
so $A\in \cA$. By~\eqref{eq:relationya} and the fact that $A\bfo = y(\vec G)$, we conclude that $y(\vec G)\in \cX$.  
\end{proof}

With Lemma~\ref{lem:necessaryyinY} above, we establish the necessity of Condition~$B'$.

\begin{proof}[Proof of necessity of Condition~$B'$ for weak $H$-property]
We show that if $x^*\notin \cX$, then~\eqref{eq:nonevent} holds. Recall that for a random graph $\vec G_n \sim W$, $x(\vec G_n) = \frac{1}{n} y(\vec G_n)$ is the empirical concentration vector of $\vec G_n$, and it converges to $x^*$~\aas.   
Since $x^* \notin \cX$ and since $\cX$ is a closed subset of $\R^m$, it holds that $x(\vec G_n)\notin \cX$ \aas. 
By Lemma~\ref{lem:necessaryyinY}, if  $x(\vec G_n)\notin \cX$ (and hence, $y(\vec G_n)\notin \cX$), then $\vec G_n$ cannot have a cycle cover. 
\end{proof}

Next, given $\vec G_n \sim W$, we define
\begin{equation}\label{eq:defomegan}
        \omega(\vec G_n):= \sqrt{n}(x(\vec G_n) - x^*) + x^*.
\end{equation}
The following result is known (see, e.g.,~\cite{belabbas2021h}):

\begin{lemma}\label{lem:limitdist}
    The random variable $\omega(\vec G_n)$ converges in distribution to the Gaussian random variable $\omega^* \sim N(x^*, \Sigma)$, 
    where $\operatorname{Diag}(x^*)$ is the diagonal matrix whose $ii$th entry is~$x^*_i$ and $\Sigma:= \operatorname{Diag}(x^*) - x^* {x^*}^\top$. 
    The rank of $\Sigma$ is $(m - 1)$ and its null space is spanned by $\bfo$.
\end{lemma}

With Lemmas~\ref{lem:necessaryyinY} and~\ref{lem:limitdist}, we establish the necessity of Condition $A$:

\begin{proof}[Proof of necessity of Condition $A$ for weak $H$-property]
We need to show that if 
\begin{equation}\label{eq:crkz1}
\crk(Z) \geq 1,
\end{equation} 
then~\eqref{eq:nonevent} holds. 
We may as well assume that Condition $B'$ holds, i.e., $x^* \in \cX$. 

To proceed, we first normalize the node-cycle incidence vectors $z_j$ so that their one-norm is~$1$: 
$$\bar z_j := \frac{z_j}{\|z_j\|_1}, \quad \mbox{for all } j = 1,\ldots, k.$$
Let $\overline\cX$ be the convex hull generated by $\bar z_1,\ldots, \bar z_k$,     
which comprises all $x\in \cX$ such that $\|x\|_1 = 1$. 
Since $\|x^*\|_1 = 1$ and since $x^*\in \cX$ (by assumption), we have that 
\begin{equation}\label{eq:xstarinbarX}
x^* \in \overline \cX. 
\end{equation}
Similarly, since $\|x(\vec G_n)\|_1 = 1$, we have that 
$x(\vec G_n) \in \cX$ if and only if $x(\vec G_n)\in \overline \cX$.

Next, let $L$  be the affine hyperplane in $\R^m$ spanned by $e_1,\ldots, e_m$, which contains the standard simplex. Let $L'$ be the affine space spanned by $\bar z_1,\ldots, \bar z_k$, which is the affine space of least dimension that contains $\overline \cX$.  
By our hypothesis~\eqref{eq:crkz1}, 
$$\dim L' \leq m - 2 < m - 1 = \dim L,$$
i.e., $L'$ is a {\it proper} affine subspace of $L$. 

We now establish a sequence of inequalities that bound from above the probability that $\vec G_n\sim W$ has a cycle cover. 
By Lemma~\ref{lem:necessaryyinY}, it is necessary that $x(\vec G_n)\in \cX$ for $\vec G_n$ to have a cycle cover, so
\begin{multline}\label{eq:onedirectionforx}
\bP(\vec G_n \sim W \mbox{ has a cycle cover}) \leq \bP(x(\vec G_n)\in \cX) \\ = \bP(x(\vec G_n) \in \overline \cX) \leq \bP(x(\vec G_n) \in L'). 
\end{multline}
Then, by~\eqref{eq:defomegan} and~\eqref{eq:xstarinbarX}, we have that 
$x(\vec G_n) \in L'$ if and only if $\omega(\vec G_n)\in L'$, so 
\begin{equation}\label{eq:twoequivalentevents}
\bP(x(\vec G_n) \in L') = \bP(\omega(\vec G_n)\in L'). 
\end{equation}
Combining~\eqref{eq:onedirectionforx} and~\eqref{eq:twoequivalentevents}, we have that 
\begin{equation}\label{eq:boundnoneventfromabove1}
\bP(\vec G_n \sim W \mbox{ has a cycle cover}) \leq \bP(\omega(\vec G_n)\in L').
\end{equation}
Finally, we appeal to Lemma~\ref{lem:limitdist} to obtain that 
$$
\lim_{n\to\infty}\bP(\omega(\vec G_n)\in L') = \lim_{n\to\infty} \bP(\omega^*\in L') = 0,
$$
where the last equality follows from the fact that $L'$ is a proper affine subspace of $L$ and the fact that the Gaussian random variable $\omega^*$ has the entire $L$ as its support.
\end{proof}

Finally, we establish the necessity of Condition~$C$: 

\begin{proof}[Proof of necessity of Condition~$C$ for strong $H$-property] 
Let $\vec G$ be an  $\vec S$-partite graph such that $y(\vec G)\in \N^m$, so  $\pi^{-1}(u_i)$ contains at least one node. If $\vec G$ has a Hamilton cycle $\vec H$, then $\pi(\vec H)$ is a closed walk of $\vec S$ that visits every node at least once, which implies that $\vec S$ is strongly connected. In other words, we have just shown that if $\vec S$ is not strongly connected and if $y(\vec G) \in \N^m$, then $\vec G$ does not have a Hamilton cycle. 
Now, let $\vec G_n\sim W$. Since $x(\vec G_n)$ converges to $x^*$ \aas~and since all the entries $x^*_i$ are positive, we have that 
$y_i(\vec G_n) = |\pi^{-1}(u_i)| = \Theta(n) \mbox{ \aas}$. 
The above arguments then imply that if $\vec S$ is not strongly connected, then \aas~$\vec G_n$ does not have a Hamilton cycle, i.e.,~\eqref{eq:posnonevent} holds.  
\end{proof}

\section{Pre-processing: Removal of self-loops}\label{sec:loopfree}

Let $W$ be a step-graphon and $\sigma = (\sigma_0, \ldots, \sigma_{m-1}, \sigma_*)$ be a partition for $W$, with $\sigma_0 = 0$ and $\sigma_* = 1$. 
Let $\vec S$ be the associated skeleton graph, and $\vec S_1,\ldots, \vec S_q$ be the SCCs of $\vec S$. 
We assume that the skeleton graph $\vec S$ associated with $W$ has a self-loop, say, on node $u_m\in V(\vec S_q)$ and that the partition $\sigma$ is fine enough such that $\vec S_q$ has at least two nodes. 

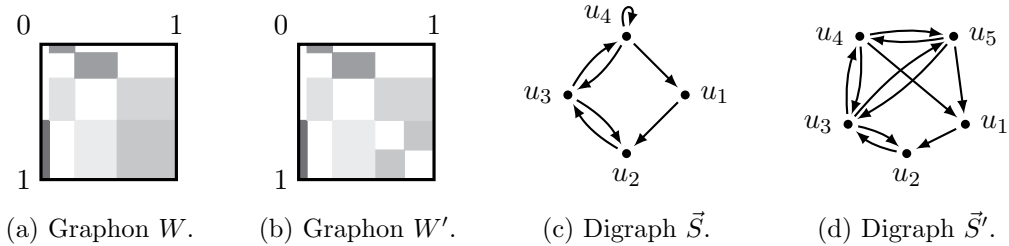
\begin{figure}[t]
\centering
    \begin{subfigure}{.22\textwidth}
\centering
\begin{tikzpicture}[scale=.45]
\filldraw [fill=Gray!30!black!80, draw=Gray!30!black!80] (0,0) rectangle (1/4,7/4); % 41
\filldraw [fill=Gray!70!black!20, draw=Gray!70!black!20] (1/4,7/4) rectangle (4/4,12/4); % 32
\filldraw [fill=Gray!70!black!70, draw=Gray!70!black!70] (4/4,12/4) rectangle (9/4,15/4); % 23
\filldraw [fill=Gray!10!black!20, draw=Gray!10!black!20] (9/4,7/4) rectangle (16/4,12/4); % 34
\filldraw [fill=Gray!60!black!70, draw=Gray!60!black!70] (1/4,15/4) rectangle (4/4,16/4); % 12
\filldraw [fill=Gray!30!black!10, draw=Gray!30!black!10] (4/4,0) rectangle (9/4,7/4); % 43
\filldraw [fill=Gray!30!black!30, draw=Gray!30!black!30] (9/4,0) rectangle (16/4,7/4); % 44

\draw [draw=black,very thick] (0,0) rectangle (4,4);

\node [above left] at (0,4) {$0$};
\node [left] at (0,0) {$1$};
\node [above] at (4,4) {$1$};
\end{tikzpicture}
\caption{Graphon $W$.}\label{sfig1:noloop}
\end{subfigure}
\begin{subfigure}{.22\textwidth}
\centering
 \begin{tikzpicture}[scale=.45]
\filldraw [fill=Gray!30!black!80, draw=Gray!30!black!80] (0,0) rectangle (1/4,7/4); % 41
\filldraw [fill=Gray!70!black!20, draw=Gray!70!black!20] (1/4,7/4) rectangle (4/4,12/4); % 32
\filldraw [fill=Gray!70!black!70, draw=Gray!70!black!70] (4/4,12/4) rectangle (9/4,15/4); % 23
\filldraw [fill=Gray!30!black!30, draw=Gray!30!black!30] (9/4,0) rectangle (25/8,7/8); % 44
\filldraw [fill=Gray!30!black!30, draw=Gray!30!black!30] (25/8,7/8) rectangle (4,7/4); % 55
\filldraw [fill=Gray!10!black!20, draw=Gray!10!black!20] (9/4,7/4) rectangle (16/4,12/4); % 34
\filldraw [fill=Gray!60!black!70, draw=Gray!60!black!70] (1/4,15/4) rectangle (4/4,16/4); % 12
\filldraw [fill=Gray!30!black!10, draw=Gray!30!black!10] (4/4,0) rectangle (9/4,7/4); % 43

\draw [draw=black,very thick] (0,0) rectangle (4,4);

\node [above left] at (0,4) {$0$};
\node [left] at (0,0) {$1$};
\node [above] at (4,4) {$1$};
\end{tikzpicture}
\caption{Graphon $W'$.}\label{sfig2:noloop}
\end{subfigure}
\begin{subfigure}{.24\textwidth}
\centering
    \begin{tikzpicture}[scale=0.78]
    \tikzset{every loop/.style={}}
		\node [circle,fill=black,inner sep=1.2pt,label=above left:{$u_4$}] (4) at (0, 0) {};
		\node [circle,fill=black,inner sep=1.2pt,label=left:{$u_3$}] (3) at (-1, -1) {};
		\node [circle,fill=black,inner sep=1.2pt,label=below:{$u_2$}] (2) at (0, -2) {};
		\node [circle,fill=black,inner sep=1.2pt,label=right:{$u_1$}] (1) at (1, -1) {};

  \path[draw,thick,shorten >=2pt,shorten <=2pt, -latex]
		 (1) edge[-latex, black] (2) % 1 2
		 (2) edge[-latex, bend left=15, black] (3) % 2 3
          (4) edge[-latex, black] (1)
          (4) edge[-latex, bend left=15, black] (3)
          (3) edge[-latex, bend left=15, black] (4)
		 (3) edge[-latex, bend left=15, black] (2) % 3 4
		 (4) edge[loop above, -latex, black] (4); % 5 5
\end{tikzpicture}
\caption{Digraph $\vec S$.}\label{sfig1:ss'}
\end{subfigure} 
\begin{subfigure}{.24\textwidth}
\centering
    \begin{tikzpicture}[scale=0.78]
		\node [circle,fill=black,inner sep=1.2pt,label=left:{$u_4$}] (4) at (-0.8, 0) {};
            \node [circle,fill=black,inner sep=1.2pt,label=right:{$u_5$}] (5) at (0.8, 0) {};
		\node [circle,fill=black,inner sep=1.2pt,label=left:{$u_3$}] (3) at (-1.0, -1.5) {};
		\node [circle,fill=black,inner sep=1.2pt,label=below:{$u_2$}] (2) at (0, -2) {};
		\node [circle,fill=black,inner sep=1.2pt,label=right:{$u_1$}] (1) at (1.0, -1.5) {};

  \path[draw,thick,shorten >=2pt,shorten <=2pt, -latex]
		 (1) edge[-latex, black, opacity = 1.0] (2) % 1 2
		 (2) edge[-latex, bend left=15, black, opacity = 1.0] (3) % 2 3
          (3) edge[-latex, bend left=15, black, opacity = 1.0] (2)
          (4) edge[-latex, black, opacity = 1.0] (1)
          (4) edge[-latex, bend left=12, black, opacity = 1.0] (3)
          (3) edge[-latex, bend left=12, black, opacity = 1.0] (4)
		 (4) edge[-latex, bend left=10, black, opacity = 1.0] (5)
          (5) edge[-latex, bend left=10, black, opacity = 1.0] (4)
          (5) edge[-latex, bend left=10, black, opacity = 1.0] (3)
          (3) edge[-latex, bend left=10, black, opacity = 1.0] (5)
          (5) edge[-latex, black, opacity = 1.0] (1) ; % 5 5
\end{tikzpicture}
\caption{Digraph $\vec S'$.}\label{sfig2:ss'}
\end{subfigure}

\caption{The step-graphon $W'$ in~(b)  is obtained from $W$ in~(a) by first subdividing the right-bottom block into $2$-by-$2$ sub-blocks of equal size and then setting the value of the two diagonal sub-blocks to zero. Let $\sigma  := \frac{1}{16}(0,1,4,9,16)$ be a partition for~$W$. The subdivision gives rise to  the partition  $\sigma' := \frac{1}{16}(0,1,4,9,12.5,16)$ for $W'$. 
The two skeleton graphs $\vec S$ and $\vec S'$, shown in~(c) and~(d), are associated with $(W,\sigma)$ and $(W',\sigma')$, respectively. 
}\label{fig:noloop}
\end{figure}

\xc{Surgery on the nonzero diagonal block of $W$:} We introduce a new step-graphon $W'$ as follows.   
Let $\sigma_m:= \frac{1}{2}(\sigma_{m-1} + 1)$ and  
\begin{equation}\label{eq:defW'}
W'(s,t) := 
\begin{cases}
0 &  \mbox{if } (s,t)\in [\sigma_{m-1}, \sigma_m)^2 \cup [\sigma_m,1]^2,\\
W(s,t) & \mbox{otherwise}. 
\end{cases}
\end{equation}
In words, $W'$ is obtained from $W$ by first subdividing the block $R_{mm} = [\sigma_{m-1},1]^2$ into four sub-blocks: 
\begin{align*}
R_{mm, 11} &:= [\sigma_{m-1},\sigma_m)^2, \quad &  R_{mm,12} & := [\sigma_{m-1},\sigma_m) \times [\sigma_m, 1], \\
R_{mm,21}&:= [\sigma_m,1] \times [\sigma_{m-1},\sigma_m),  & R_{mm, 22} &:= [\sigma_m, 1]^2.
\end{align*}
and then, setting the value of $W(s,t)$ to $0$ if $(s,t)\in R_{mm,11} \cup R_{mm,22}$ while keeping $W(s,t)$ unchanged otherwise. See Figure~\ref{fig:noloop} for illustration.

The goal of this section is to show that $W'$ inherits any Condition~$\star$, for $\star \in \{A, B, C\}$, that is satisfied by $W$. Precisely, we have

\begin{theorem}\label{thm:reducetozerodiagonal}
Let $W$ and $W'$ be given as above. If $W$ satisfies Condition~$\star$, for $\star \in\{ A, B, C\}$, then $W'$ satisfies the same condition.
\end{theorem}

Let $\sigma':= (\sigma_0,\ldots,\sigma_{m-1}, \sigma_m, \sigma_*)$. It is clear that $\sigma'$ is a partition for $W'$.  
Let $x'^*$, $\vec S'$, $Z'$, and $\cX'$ be the concentration vector, the skeleton graph, the node-cycle incidence matrix, and the node-circulation cone of $W'$ for $\sigma'$, respectively. With slight abuse of terminology, we say that $\vec S'$ is obtained from $\vec S$ by performing the surgery on node $u_m$.   
It is clear that $\vec S'$ has one less self-loop than $\vec S$ does. 

If $W'$ still has a nonzero diagonal block, then we perform the surgery again for $W'$ on the corresponding block. Iterating this procedure until we obtain a graphon whose diagonal blocks are zero. We introduce the following definition: 

\begin{definition}\label{def:loopfree}
A step-graphon $W$ is {\bf loop free} if there is a (and hence, any) partition such that the associated skeleton graph does not have any self-loop. 
\end{definition}
 
The following result is then an immediate consequence of Theorem~\ref{thm:reducetozerodiagonal}: 

\begin{corollary}\label{cor:loopfreegraphon}
If a step-graphon $W$ satisfies Condition~$\star$, for $\star \in \{A, B, C\}$, then there is a loop-free step-graphon $W'$ such that $W'\leq W$ and satisfies the same condition. 
\end{corollary}

The remainder of the section is devoted to the proof of Theorem~\ref{thm:reducetozerodiagonal}.  
We deal with the three conditions in the order of $C$, $A$, and $B$ in three subsections.

\subsection{Proof for Condition~$C$}
In this subsection, we show that
$$
\vec S \mbox{ is strongly connected (with at least 2 nodes)} \quad \Longrightarrow \quad \vec S' \mbox{ is strongly connected}.
$$

First, note that by~\eqref{eq:defW'}, the digraph $\vec S'$ can be obtained from $\vec S$ by first adding a new node $u_{m+1}$ and the following set of new edges: 
\begin{equation}\label{eq:defe0fors'}
E_{m+1} := \{u_iu_{m+1} \mid u_iu_m \in E(\vec S)\} \cup \{u_{m+1}u_j \mid  u_mu_j \in E(\vec S)\},  
\end{equation}
and then deleting the self-loop $u_m u_m$.  
Precisely, 
\begin{equation}\label{eq:vefors'}
V(\vec S') = V(\vec S) \cup \{u_{m+1}\} \quad \mbox{and} \quad E(\vec S') = E(\vec S) \cup E_{m+1} - \{u_m u_m\}.
\end{equation}

We now show that for any two distinct nodes $u_i, u_j\in V(\vec S')$, there is a walk from $u_i$ to $u_j$. Consider the following three cases:

{\it Case 1: $u_i\neq u_{m+1}$ and $u_j\neq u_{m+1}$.} Since $\vec S$ is strongly connected, there is a path $\vec P$ from $u_i$ to $u_j$ in $\vec S$. By~\eqref{eq:defe0fors'} and~\eqref{eq:vefors'}, $\vec P$ is also a path of $\vec S'$.

{\it Case 2: $u_i = u_{m+1}$.} Let $u_{i_1}\cdots u_{i_\ell}$, with $u_{i_1} = u_m$ and $u_{i_\ell} = u_j$, be a walk of $\vec S$ from $u_m$ to $u_j$. In the case $u_j = u_m$, the walk is closed (such a closed walk exists because $\vec S$ is strongly connected and has $m \geq 2$ nodes).   
Since $u_mu_{i_2}\in E(\vec S)$, by~\eqref{eq:defe0fors'} $u_{m+1}u_{i_2}\in E(\vec S')$ and hence, $u_{m+1}u_{i_2}\cdots u_{i_\ell}$ is a walk of $\vec S'$ from $u_{m+1}$ to $u_j$. 

{\it Case 3: $u_j = u_{m+1}$.} Similarly, if $u_{j_1}\cdots u_{j_\ell}$ is a walk of $\vec S$, with $u_{j_1} = u_i$ and $u_{j_\ell} = u_m$,  then $u_{j_1}\cdots u_{j_{\ell-1}}u_{m+1}$ is a walk of $\vec S'$ from $u_i$ to $u_{m+1}$. \hfill{\qed}

\subsection{Proof for Condition $A$}
In this subsection, we show that
$$
\crk(Z) = 0 \quad \Longrightarrow \quad \crk(Z') = 0.
$$

Recall that $\vec S_1,\ldots, \vec S_q$ are the SCCs of $\vec S$. Let 
$\vec S'_p := \vec S_p$, for $p = 1,\ldots, q-1$, and $\vec S'_q$ be obtained from $\vec S_q$ by performing the surgery on the node $u_m$. It follows from~\eqref{eq:defe0fors'} and~\eqref{eq:vefors'} that $\vec S'_1,\ldots, \vec S'_q$ are the SCCs of $\vec S'$.

Let $B_{\vec S'}$ be the bipartite graph associated with $\vec S'$, whose node set is a disjoint union of
$V'(B_{\vec S'}) = V'(B_{\vec S}) \cup \{u'_{m+1}\}$ and 
$V''(B_{\vec S'}) = V''(B_{\vec S}) \cup \{u''_{m+1}\}$, and whose edge set is given by
\begin{multline}\label{eq:defebs}
E(B_{\vec S'}) = E(B_{\vec S}) \cup \{(u'_i, u''_{m+1}) \mid (u'_i, u''_m)\in E(B_{\vec S}) \} \\ 
\cup  \{(u'_{m+1}, u''_i)\mid (u'_m, u''_i)\in E(B_{\vec S}) \}
- \{(u'_m,u''_m)\}.
\end{multline}

Since $\crk(Z) = 0$, it follows from Lemma~\ref{lem:connectednessofB} that every  bipartite graph $B_{\vec S_p}$, for $p = 1,\ldots, q$, is connected. Using the same lemma, we have that $\crk(Z') = 0$ if and only if $B_{\vec S'_p}$ is connected for all $p = 1,\ldots, q$. Since $\vec S'_p = \vec S_p$ for $p = 1,\ldots, q - 1$, it suffices to show that $B_{\vec S'_q}$ is connected. 

By the above arguments, we can assume without loss of generality that $\vec S$ is itself strongly connected (and $\vec S'$ is obtained from $\vec S$ by performing the surgery on $u_m$). Under the assumption, we establish the following lemma: 

\begin{lemma}
If $B_{\vec S}$ is connected, then so is $B_{\vec S'}$.  
\end{lemma}

\begin{proof}
We show that for any $u'_i\in V'(B_{\vec S'})$ and any $u''_j\in V''(B_{\vec S'})$, there is a path of $B_{\vec S'}$ from $u'_i$ to $u''_j$. Consider the following four cases: 

{\it Case 1: $u'_i \neq u'_{m+1}$ and $u''_j \neq u''_{m+1}$.} 
Let $P$ be a path of $B_{\vec S}$ that connects $u'_i$ and $u''_j$. If the path does not contain the edge $(u'_m, u''_m)$, then $P$ is also a path of $B_{\vec S'}$. We thus assume that $P$ contains $(u'_m, u''_m)$. 
Since $m \geq 2$ and since $B_{\vec S}$ is connected, at least one of the two nodes $u'_m$ and $u''_m$ has degree at least $2$ within $B_{\vec S}$. 
Without loss of generality, we assume that $\deg(u'_m) \geq 2$ and that $(u'_m, u''_\ell)$, with $u''_\ell \neq u''_m$, is an edge of $B_{\vec S}$. 
By~\eqref{eq:defebs}, we have that $(u'_m, u''_\ell)$, $(u'_{m+1}, u''_\ell)$, and $(u'_{m+1}, u''_m)$ are edges of $B_{\vec S'}$. 
Replacing the segment $u'_mu''_m$ in $P$ with $u'_mu''_\ell u'_{m+1} u''_m$, we obtain a walk of $ B_{\vec S'}$ that connects $u'_i$ and $u''_j$. 

{\it Case 2: $u'_i = u'_{m+1}$ and $u''_j = u''_{m+1}$.} By the same arguments in Case 1, we can assume without loss of generality that $(u'_m, u''_\ell)$, with $u''_\ell \neq u''_m$, is an edge of $B_{\vec S}$. Then, $u'_{m+1}u''_\ell u'_m u''_{m+1}$ is a path from $u'_{m+1}$ to $u''_{m+1}$. 

{\it Case 3: $u'_i \neq u'_{m+1}$ and $u''_j = u''_{m+1}$.} Let $P$ be a path of $B_{\vec S}$ from $u'_i$ to $u''_m$. Replacing the last node $u''_m$ of $P$ with $u''_{m+1}$, we obtain a path of $B_{\vec S'}$ from $u'_i$ to $u''_{m+1}$. 

{\it Case 4: $u'_i = u'_{m+1}$ and $u''_j \neq u''_{m+1}$.} Similarly, let $P$ be a path of $B_{\vec S}$ from $u'_m$ to $u''_j$. Replacing the first node $u'_m$ of $P$ with $u'_{m+1}$, we obtain a path of $B_{\vec S'}$ from $u'_{m+1}$ to $u''_j$. 
\end{proof}

\subsection{Proof for Condition $B$}\label{ssec:proofthatx'inY'noloop}
In this subsection, we show that
$$
x^*\in \rint \cX \quad \Longrightarrow \quad x'^* \in \rint \cX'.
$$

We start by relating the cycles of $\vec S'$ to those of $\vec S$. 
Label the cycles of $\vec S$ in a way such that the first $\ell$ cycles 
$\vec C_1,\ldots, \vec C_\ell$, for some $\ell \leq k$, contain the node $u_m$ and that $\vec C_1 = u_mu_m$ is the self-loop. 

The self-loop $\vec C_1$ induces the $2$-cycle 
$\vec C'_1 := u_m u_{m+1} u_m$ of $\vec S'$.  
Each cycle $\vec C_p$, for $2 \leq p \leq \ell$, induces four different cycles of $\vec S'$ as follows: $\vec C_{p, 1} := \vec C_p$ and $\vec C_{p, 2}$, $\vec C_{p, 3}$, $\vec C_{p, 4}$ are obtained from $\vec C_p$ by substituting the node $u_m$ with $u_{m+1}$, $u_mu_{m+1}$, $u_{m+1}u_m$, respectively. 
Thus, the set of cycles of $\vec S'$ is given by 
$$\{\vec C'_1\} \cup \{\vec C_{p, i} \mid 2\leq p\leq \ell \mbox{ and } 1\leq i \leq 4\} \cup \{\vec C_q \mid  \ell + 1 \leq q \leq k\}.$$

To illustrate, consider the digraph $\vec S$ in Figure~\ref{sfig1:ss'} and the corresponding digraph $\vec S'$ in Figure~\ref{sfig2:ss'}. 
The digraph $\vec S$ has $4$ cycles as exhibited in~\eqref{eq:cyclesofexampleS}.  
The first three cycles contain the node $u_4$. 
The self-loop $\vec C_1$ induces the $2$-cycle $\vec C'_1 = u_4u_5u_4$ in $\vec S'$. The cycle $\vec C_2$ induces four cycles of $\vec S'$, which are 
$\vec C_{2,1} = u_3u_4u_3$, $\vec C_{2,2} = u_3u_5u_3$, $\vec C_{2,3} = u_3u_4u_5u_3$, and $\vec C_{2,4} = u_3u_5u_4u_3$. 
Similarly, the four cycles of $\vec S'$ induced by $\vec C_3$ are
    $\vec C_{3,1}  = u_1u_2u_3u_4u_1$,  $\vec C_{3,2}  = u_1u_2u_3u_5u_1$, 
    $\vec C_{3,3} = u_1u_2u_3u_4u_5u_1$,  and $\vec C_{3,4}  = u_1u_2u_3u_5u_4u_1$. 
Thus, $\vec S'$ has ten cycles $\vec C'_1$, $\vec C_{2,1},\ldots, \vec C_{2,4}$, $\vec C_{3,1},\ldots, \vec C_{3,4}$, and $\vec C_4$.

Let $z'_1$, $z'_{p, i}$, and $z'_q$ be the node-cycle incidence vectors of $\vec S'$ corresponding to $\vec C'_1$, $\vec C_{p, i}$, and $\vec C_q$, respectively.  
To relate these vectors to the $z_j$'s, we first augment each $z_j$ by adding a zero entry at the end. Precisely, we define 
$
\hat z_j :=
(z_j; 0)
 \in \R^{m+1}$, for all $j = 1,\ldots, k$. 
Then, 
\begin{equation}\label{eq:relatez'toz}
\left\{
\begin{aligned}
& z'_1  = e_m + e_{m+1}, & \\
& z'_{p, 1} = \hat z_p, \, z'_{p, 2} = \hat z_p - e_m + e_{m+1}, \, z'_{p, 3} = z'_{p, 4} = \hat z_p + e_{m+1},  & \mbox{for } 2\leq p\leq \ell, \\
& z'_q = \hat z_q,  & \mbox{for } \ell < q \leq k.
\end{aligned}
\right.
\end{equation}

Note that
$z'_{p, 3} = z'_{p, 4} = \frac{1}{2} (z'_1 + z'_{p, 1} + z'_{p, 2})$, 
which implies that $z'_{p, 3}$ and $z'_{p'_4}$ are {\em not} extremal generators of $\cX'$. 
Thus, it suffices to show that $x'^*$ can be expressed as a {\em positive} combination of $z'_1$, $z'_{p, 1}$, $z'_{p,2}$, and $z'_q$.   
First, by the definition of $\sigma'$, we have that
$x'^* = (x^*_1, \cdots, x^*_{m-1}, \frac{x^*_m}{2}, \frac{x^*_m}{2})$. 
Let  $\hat x^* := (x^*; 0)$; then, we can express $x'^*$ as
\begin{equation}\label{eq:defxprimestar}
x'^* = \hat x^* + \frac{x^*_m}{2} (e_{m + 1} - e_m).
\end{equation}
Since $x^*\in \rint \cX$, there exist positive coefficients $c_j$'s such that $x^* = \sum_{j = 1}^k c_j z_j$. It follows from the definitions of $\hat z_j$ and of $\hat x^*$ that
\begin{equation}\label{eq:xhatzhat}
\hat x^* = \sum_{j = 1}^k c_j \hat z_j. 
\end{equation}
Since $\vec C_1,\ldots, \vec C_\ell$ are the cycles of $\vec S$ that contain $u_m$,
\begin{equation}\label{eq:sumforxm}
x^*_m = \sum_{j = 1}^\ell c_j.
\end{equation}
We define positive coefficients as follows: 
\begin{equation}\label{eq:defc'}
\begin{cases}
c'_1  := c_1/2, & \\
c'_{p, i}  :=  c_p/2 & \quad  \mbox{for } 2\leq p \leq \ell \mbox{ and } 1 \leq i \leq 2, \\ 
c'_q := c_q  & \quad \mbox{for } \ell + 1 \leq q \leq k.
\end{cases}
\end{equation}
Then, following~\eqref{eq:defxprimestar}, we have that 
$$
\begin{aligned}
x'^* & = \hat x^* + \frac{x^*_m}{2} (e_{m + 1} - e_m) \\
& = \sum_{j = 1}^k c_j \hat z_j + \frac{x^*_m}{2} (e_{m + 1} - e_m) \\
& = \sum_{q = \ell + 1}^k c'_q \hat z_q  +  \sum_{p = 2}^\ell \left [\sum_{i = 1}^2 c'_{p, i} \right ] \hat z_p   + c_1 e_m  + \frac{x^*_m}{2} (e_{m + 1} - e_m) \\
& = \sum_{q = \ell + 1}^k c'_q z'_q  +  \sum_{p = 2}^\ell \sum_{i = 1}^2 c'_{p, i} z'_{p, i} + c_1 e_m  + \frac{1}{2}\left [x^*_m - \sum_{j = 2}^\ell c_j \right ](e_{m+1} - e_m) \\
& = \sum_{q = \ell + 1}^k c'_q z'_q  +  \sum_{p = 2}^\ell \sum_{i = 1}^2 c'_{p, i} z'_{p, i} + \frac{1}{2} c_1 (e_{m + 1} + e_m) \\
& =  \sum_{q = \ell + 1}^k c'_q z'_q +  \sum_{p = 2}^\ell \sum_{i = 1}^2 c'_{p, i} z'_{p, i}   + c'_1 z'_1,  
\end{aligned}
$$
where the second equality follows from~\eqref{eq:xhatzhat}, the third equality follows from~\eqref{eq:defc'}, the fourth equality follows from~~\eqref{eq:relatez'toz} and~\eqref{eq:defc'}, the fifth equality follows from~\eqref{eq:sumforxm}, and the last equality follows from~\eqref{eq:relatez'toz} and~\eqref{eq:defc'}. This completes the proof. \hfill{\qed}

\section{Hamiltonicity of complete $\vec S$-partite graphs}\label{sec:constructHam}

In this section, we assume that $\vec S$ does {\it not} have a self-loop and that Condition $B$ ($x^*\in \rint \cX$) is satisfied.  
Let $\cU$ be an open neighborhood of $x^*$ in $\cX$. Then, there is a continuous function 
$\gamma: \cU \to \R_{>0}^k$ such that 
\begin{equation}\label{eq:defphi}
x = Z \gamma(x) \quad \mbox{for all } x\in \cU. 
\end{equation}
Let 
$\gamma_0:= \frac{1}{2}\min_{j = 1}^k \gamma_j(x^*)$. 
Shrink $\cU$ if necessary so that 
\begin{equation}\label{eq:defepsilon}
\gamma_j(y)  >  \gamma_0, \quad \mbox{for all } j= 1,\ldots, k \mbox{ and for all } y \in \cU.
\end{equation}
We then introduce the following subset of $\cX$:
\begin{equation}\label{eq:defY0}
\cX_0 := \{y\in \cX \cap \N^m\mid \|y\|_1 \geq 1/\gamma_0   \quad \mbox{and} \quad y/\|y\|_1 \in \cU \}.
\end{equation}
In words, $\cX_0$ is the set of all integer-valued $y\in \cX$ such that $\|y\|_1$ is sufficiently large and $y/\|y\|_1$ is sufficiently close to $x^*$. 

The main result of this section is a sufficient condition for a complete $\vec S$-partite graph to have a cycle cover or a Hamilton cycle. We state it below:

\begin{theorem}\label{thm:constructHam}
    The following two items hold: 
    \begin{enumerate}
    \item For any integer-valued $y\in \cX$, $\vec K_y$ has a cycle cover. 
    \item If $\vec S$ is strongly connected, then for any $y\in \cX_0$, $\vec K_y$ has a Hamilton cycle. 
    \end{enumerate} 
\end{theorem}

By Lemma~\ref{lem:necessaryyinY}, if $\vec K_y$ has a cycle cover, then $y\in \cX$. Combining this fact with item~1 of the above theorem, we have that $y\in \cX \cap \N_0^m$ is both necessary and sufficient for $\vec K_y$ to have a cycle cover. 

We establish the two items of Theorem~\ref{thm:constructHam} in two subsections. 

\subsection{Proof of item~1 of Theorem~\ref{thm:constructHam}}
We start by decomposing $y\in \cX$ into an integer combination of the node-cycle incidence vectors $z_j$. This is feasible as we show in the following lemma:

\begin{lemma}\label{lem:existenceofHD1}
For any integer-valued $y\in \cX$, there exist $c_1,\ldots, c_k\in \N_0$ such that 
    \begin{equation}\label{eq:step1}
    y = \sum_{j = 1}^k c_j z_j.
    \end{equation}
\end{lemma}

\begin{proof}
Since $y$ is integer valued, it follows from the integrality theorem that there exists an integer-valued $A\in \cA$ such that 
\begin{equation}\label{eq:y''toA''}
y = A \bfo.
\end{equation} 
We show that there exist $c_1,\ldots, c_k\in \N_0$ such that 
\begin{equation}\label{eq:decomposeA}
    A = \sum_{j = 1}^k c_j A_j.
\end{equation}
Since $A\in \cA$, there exist $r_1,\ldots,r_{k}\in \R_{\geq 0}$ such that
$A = \sum_{j = 1}^k r_j A_j$.  
Since $A$ is integer valued, it holds that if $r_j > 0$ for some $j= 1,\ldots, k$, then $A' := (A - A_j)$ has nonnegative entries and is integer valued. We claim that $A'\in \cA$. To wit, note that $\supp(A')\subseteq \supp(A)$ and $\supp(A)\subseteq E(\vec S)$, so $\supp(A')\subseteq E(\vec S)$. 
Also, note that 
$$A' \bfo = A\bfo - A_j \bfo = A^\top \bfo - A_j^\top \bfo = A'^\top \bfo.$$
This establishes the claim. If $A' \neq 0$, then we can repeat the same arguments to find some $j'= 1,\ldots, k$ such that $(A' - A_{j'})\in \cA$. This iteration will terminate in finite steps and we obtain~\eqref{eq:decomposeA}. 
Now, using~\eqref{eq:y''toA''},~\eqref{eq:decomposeA}, and the fact that $A_j \bfo = z_j$, we obtain~\eqref{eq:step1}.
\end{proof}

Using the decomposition~\eqref{eq:step1}, we exhibit a desired cycle cover of $\vec K_y$ in the following lemma:

\begin{lemma}\label{lem:existenceofHD2}
    Let $c_1,\ldots, c_k\in \N_0$ be given as in Lemma~\ref{lem:existenceofHD1} so that~\eqref{eq:step1} holds. Then, $\vec K_y$ has a cycle cover $\vec H$, which   contains, for each $j = 1,\ldots, k$, $c_j$ cycles that are isomorphic to $\vec C_j$ under the map $\pi$.  
\end{lemma}

\begin{proof}
The proof will be carried out by induction on $c:=\sum_{j = 1}^k c_j$. 
For the base case $c = 0$, item~2 holds trivially. 
For the inductive step, we assume that item~2 holds for $(c - 1) \geq 0$ and prove for $c$. 

Without loss of generality, we assume that $c_1 \geq 1$ and write $\vec C_1 = u_1 u_2 \cdots u_{d_1} u_1$, where $d_1$ is the length of $\vec C_1$. Since $\vec S$ does not have a self-loop, $d_1 \geq 2$.  
It follows from~\eqref{eq:step1} that for each $i = 1,\ldots, d_1$, $y_i = |\pi^{-1}(u_i)| \geq 1$, so there is at least a node, say $v_i$, contained in $\pi^{-1}(u_i)$.   
Because $\vec K_y$ is complete $\vec S$-partite and because $\vec C_1$ is a cycle of $\vec S$, 
we have that
$\vec D_1:= v_1v_2\cdots v_{d_1}v_1$ is a cycle of $\vec K_y$. 
It is clear that $\vec D_1$ is isomorphic to $\vec C_1$ under the map $\pi$. 

We now remove $\vec D_1$ from $\vec K_y$ and the edges incident to $\vec D_1$. Then, the resulting graph is the complete $\vec S$-partite graph $\vec K_{y'}$, where 
$$y':= y - z_1 = (c_1 - 1) z_1 + \sum_{j = 2}^k c_j  z_j.$$
By the induction hypothesis, $\vec K_{y'}$ has a cycle cover $\vec H'$ which contains $(c_1 - 1)$ cycles isomorphic to $\vec C_1$ and $c_j$ cycles isomorphic to $\vec C_j$ for $j = 2,\ldots, k$, under the map $\pi$. 
Taking the union of $\vec H'$ and the cycle $\vec D_1$, we obtain the desired cycle cover for $\vec K_y$.
\end{proof}

\subsection{Proof of item~2 of Theorem~\ref{thm:constructHam}}
Under the assumption that $\vec S$ is strongly connected and $y\in \cX_0$, 
the two lemmas we established in the previous subsection can be strengthened. 
We first have the following result, as a strengthened version of Lemma~\ref{lem:existenceofHD1}:

\begin{lemma}
For any $y\in \cX_0$, there exist positive integers $c_1,\ldots, c_k$ such that
\begin{equation}\label{eq:step11}
    y = \sum_{j = 1}^k c_j z_j.
\end{equation}
\end{lemma}

\begin{proof}
For convenience, let $n:= \|y\|_1$. Since $y\in \cX_0$, $\frac{y}{n} \in \cU$. 
By~\eqref{eq:defphi}, we write 
\begin{equation}\label{eq:defbarA}
y = \sum_{j = 1}^k n\gamma_j(x) z_j.
\end{equation}   
Let 
$c'_j  := \left \lfloor n \gamma_j(x) \right \rfloor$ and  
$r'_j :=  n \gamma_j(x) - c'_j$,  for all $j= 1,\ldots, k$. 
By~\eqref{eq:defepsilon},~\eqref{eq:defY0}, and the hypothesis that $y\in \cX_0$, 
we have that 
$c'_j \geq \left \lfloor n\gamma_0 \right \rfloor \geq 1$, for all $j = 1,\ldots, k$. 
If $r'_j = 0$ for all $j$, then we set $c_j:= c'_j$ and~\eqref{eq:step11} holds.  Otherwise, let
\begin{equation}\label{eq:defy'y''}
y':= \sum_{j = 1}^k c'_j  z_j \quad \mbox{and} \quad y'':= y - y' = \sum_{j = 1}^k r'_j z_j. 
\end{equation}
It is clear that both $y'$ and $y''$ are integer valued and belong to $\cX$. By Lemma~\ref{lem:existenceofHD1}, there exist $c''_1,\ldots,c''_k\in \N_0$ such that 
\begin{equation}\label{eq:expressy''}
y'' = \sum_{j = 1}^k c''_j z_j.
\end{equation}
We then let $c_j:= c'_j + c''_j$ for $j= 1,\ldots, k$. It is clear that all the $c_j$'s are positive.   
Using~\eqref{eq:defy'y''} and~\eqref{eq:expressy''}, we conclude that~\eqref{eq:step11} holds.
\end{proof}

We now show that whenever $\vec S$ is strongly connected and $y$ can be expressed as~\eqref{eq:step11}, with $c_1,\ldots, c_k$ positive integers, the digraph $\vec K_y$ has a Hamilton cycle. 

\begin{lemma}
Suppose that $\vec S$ is strongly connected and that~\eqref{eq:step11} holds for some positive integers $c_1,\ldots, c_k$; then, $\vec K_y$ has a Hamilton cycle.
\end{lemma}

\begin{proof}
The proof will be carried out by induction on the number~$k$ of cycles of $\vec S$. 

\xc{Base case $k = 1$.} In this case, $\vec S$ is itself a cycle. We write $\vec S = u_1u_2\cdots u_m u_1$, for $m\geq 2$.  
By Lemma~\ref{lem:existenceofHD2}, there exists a cycle cover $\vec H$ of $\vec K_y$ which comprises~$c_1$ cycles that are isomorphic to~$\vec S$ under~$\pi$. We label these cycles as $\vec D_1,\ldots, \vec D_{c_1}$ and write 
$$
\vec D_j = v_{j,1}v_{j,2}\cdots v_{j,m} v_{j,1}, \quad \mbox{for all } j = 1,\ldots, c_1,  
$$
where the nodes are labeled such that 
$$
\pi^{-1}(u_i) = \{v_{j,i} \mid j = 1,\ldots, c_1\}, \quad \mbox{for all } i = 1,\ldots, m.
$$
Since $\vec K_y$ is complete $\vec S$-partite, we have that 
$v_{i,m}v_{j,1}$ is an edge of $\vec K_y$ for any $1\leq i, j \leq c_1$. It follows that
$$
\vec H := v_{1,1}\cdots v_{1,m}v_{2,1}\cdots v_{2,m} v_{3,1} \cdots v_{c_1,m} v_{1,1}
$$
is a Hamilton cycle of $\vec K_y$.

\xc{Inductive step.} We assume that the lemma holds for any $k'\leq k -1 $ and prove for $k$. Since $\vec S$ is strongly connected,  $\vec S$ admits an {\it  ear decomposition}. See, e.g.,~\cite[Chapter 7.2]{bang2008digraphs} and also Figure~\ref{fig:eardec} for an illustration. 
In particular, $\vec S$ can be obtained by gluing an ear $\vec P = u_{1} \cdots u_{r}$ to a strongly connected subgraph $\vec S'$, where the starting node $u_{1}$ and the ending node $u_{r}$ of the ear are nodes of $\vec S'$ while the other nodes of the ear do not belong to $\vec S'$. Note that $u_{1}$ and $u_{r}$ can be the same (in this case, $\vec P$ is a cycle). 

\begin{figure}
\centering
\subfloat[\label{sfig0:eardec}]{
\centering
    \begin{tikzpicture}[scale=0.6]
		\node [circle,fill=black,inner sep=1.2pt,label=left:{}] (4) at (-0.8, 0) {};
		\node [circle,fill=black,inner sep=1.2pt,label=left:{}] (3) at (-1.0, -1.5) {};
		\node  () at (0, -2.2) {};

  \path[draw,thick,shorten >=2pt,shorten <=2pt, -latex]
          (4) edge[-latex, bend left=12, black, opacity = 1.0] (3)
          (3) edge[-latex, bend left=12, black, opacity = 1.0] (4);
\end{tikzpicture}
}
\subfloat[\label{sfig1:eardec}]{
\centering
    \begin{tikzpicture}[scale=0.6]
		\node [circle,fill=black,inner sep=1.2pt,label=left:{}] (4) at (-0.8, 0) {};
		\node [circle,fill=black,inner sep=1.2pt,label=left:{}] (3) at (-1.0, -1.5) {};
		\node [circle,fill=red,inner sep=1.2pt,label=below:{}] (2) at (0, -2) {};

  \path[draw,thick,shorten >=2pt,shorten <=2pt, -latex]
	  (2) edge[-latex, bend left=15, red, opacity = 1.0] (3) % 2 3
          (3) edge[-latex, bend left=15, red, opacity = 1.0] (2)
          (4) edge[-latex, bend left=12, black, opacity = 1.0] (3)
          (3) edge[-latex, bend left=12, black, opacity = 1.0] (4);
\end{tikzpicture}
}
\subfloat[\label{sfig2:eardec}]{
\centering
    \begin{tikzpicture}[scale=0.6]
		\node [circle,fill=black,inner sep=1.2pt,label=left:{}] (4) at (-0.8, 0) {};
           	\node [circle,fill=black,inner sep=1.2pt,label=left:{}] (3) at (-1.0, -1.5) {};
		\node [circle,fill=black,inner sep=1.2pt,label=below:{}] (2) at (0, -2) {};
		\node [circle,fill=red,inner sep=1.2pt,label=right:{}] (1) at (1.0, -1.5) {};

  \path[draw,thick,shorten >=2pt,shorten <=2pt, -latex]
		 (1) edge[-latex, red, opacity = 1.0] (2) % 1 2
		 (2) edge[-latex, bend left=15, black, opacity = 1.0] (3) % 2 3
          (3) edge[-latex, bend left=15, black, opacity = 1.0] (2)
          (4) edge[-latex, red, opacity = 1.0] (1)
          (4) edge[-latex, bend left=12, black, opacity = 1.0] (3)
          (3) edge[-latex, bend left=12, black, opacity = 1.0] (4);
\end{tikzpicture}
}
\subfloat[\label{sfig3:eardec}]{
\centering
    \begin{tikzpicture}[scale=0.6]
		\node [circle,fill=black,inner sep=1.2pt,label=left:{}] (4) at (-0.8, 0) {};
            \node [circle,fill=red,inner sep=1.2pt,label=right:{}] (5) at (0.8, 0) {};
		\node [circle,fill=black,inner sep=1.2pt,label=left:{}] (3) at (-1.0, -1.5) {};
		\node [circle,fill=black,inner sep=1.2pt,label=below:{}] (2) at (0, -2) {};
		\node [circle,fill=black,inner sep=1.2pt,label=right:{}] (1) at (1.0, -1.5) {};

  \path[draw,thick,shorten >=2pt,shorten <=2pt, -latex]
		 (1) edge[-latex, black, opacity = 1.0] (2) % 1 2
		 (2) edge[-latex, bend left=15, black, opacity = 1.0] (3) % 2 3
          (3) edge[-latex, bend left=15, black, opacity = 1.0] (2)
          (4) edge[-latex, black, opacity = 1.0] (1)
          (4) edge[-latex, bend left=12, black, opacity = 1.0] (3)
          (3) edge[-latex, bend left=12, black, opacity = 1.0] (4)
		 (4) edge[-latex, bend left=10, red, opacity = 1.0] (5)
          (5) edge[-latex, red, opacity = 1.0] (1) ; % 5 5
\end{tikzpicture}
}
\subfloat[\label{sfig4:eardec}]{
\centering
    \begin{tikzpicture}[scale=0.6]
		\node [circle,fill=black,inner sep=1.2pt,label=left:{}] (4) at (-0.8, 0) {};
            \node [circle,fill=black,inner sep=1.2pt,label=right:{}] (5) at (0.8, 0) {};
		\node [circle,fill=black,inner sep=1.2pt,label=left:{}] (3) at (-1.0, -1.5) {};
		\node [circle,fill=black,inner sep=1.2pt,label=below:{}] (2) at (0, -2) {};
		\node [circle,fill=black,inner sep=1.2pt,label=right:{}] (1) at (1.0, -1.5) {};

  \path[draw,thick,shorten >=2pt,shorten <=2pt, -latex]
		 (1) edge[-latex, black, opacity = 1.0] (2) % 1 2
		 (2) edge[-latex, bend left=15, black, opacity = 1.0] (3) % 2 3
          (3) edge[-latex, bend left=15, black, opacity = 1.0] (2)
          (4) edge[-latex, black, opacity = 1.0] (1)
          (4) edge[-latex, bend left=12, black, opacity = 1.0] (3)
          (3) edge[-latex, bend left=12, black, opacity = 1.0] (4)
		 (4) edge[-latex, bend left=10, black, opacity = 1.0] (5)
          (5) edge[-latex, bend left=10, red, opacity = 1.0] (3)
          (5) edge[-latex, black, opacity = 1.0] (1) ; % 5 5
\end{tikzpicture}
}
\subfloat[\label{sfig5:eardec}]{
\centering
    \begin{tikzpicture}[scale=0.6]
		\node [circle,fill=black,inner sep=1.2pt,label=left:{}] (4) at (-0.8, 0) {};
            \node [circle,fill=black,inner sep=1.2pt,label=right:{}] (5) at (0.8, 0) {};
		\node [circle,fill=black,inner sep=1.2pt,label=left:{}] (3) at (-1.0, -1.5) {};
		\node [circle,fill=black,inner sep=1.2pt,label=below:{}] (2) at (0, -2) {};
		\node [circle,fill=black,inner sep=1.2pt,label=right:{}] (1) at (1.0, -1.5) {};

  \path[draw,thick,shorten >=2pt,shorten <=2pt, -latex]
		 (1) edge[-latex, black, opacity = 1.0] (2) % 1 2
		 (2) edge[-latex, bend left=15, black, opacity = 1.0] (3) % 2 3
          (3) edge[-latex, bend left=15, black, opacity = 1.0] (2)
          (4) edge[-latex, black, opacity = 1.0] (1)
          (4) edge[-latex, bend left=12, black, opacity = 1.0] (3)
          (3) edge[-latex, bend left=12, black, opacity = 1.0] (4)
		 (4) edge[-latex, bend left=10, black, opacity = 1.0] (5)
          (5) edge[-latex, bend left=10, red, opacity = 1.0] (4)
          (5) edge[-latex, bend left=10, black, opacity = 1.0] (3)
          (5) edge[-latex, black, opacity = 1.0] (1) ; % 5 5
\end{tikzpicture}
}
\subfloat[\label{sfig6:eardec}]{
\centering
    \begin{tikzpicture}[scale=0.6]
		\node [circle,fill=black,inner sep=1.2pt,label=left:{}] (4) at (-0.8, 0) {};
            \node [circle,fill=black,inner sep=1.2pt,label=right:{}] (5) at (0.8, 0) {};
		\node [circle,fill=black,inner sep=1.2pt,label=left:{}] (3) at (-1.0, -1.5) {};
		\node [circle,fill=black,inner sep=1.2pt,label=below:{}] (2) at (0, -2) {};
		\node [circle,fill=black,inner sep=1.2pt,label=right:{}] (1) at (1.0, -1.5) {};

  \path[draw,thick,shorten >=2pt,shorten <=2pt, -latex]
		 (1) edge[-latex, black, opacity = 1.0] (2) % 1 2
		 (2) edge[-latex, bend left=15, black, opacity = 1.0] (3) % 2 3
          (3) edge[-latex, bend left=15, black, opacity = 1.0] (2)
          (4) edge[-latex, black, opacity = 1.0] (1)
          (4) edge[-latex, bend left=12, black, opacity = 1.0] (3)
          (3) edge[-latex, bend left=12, black, opacity = 1.0] (4)
		 (4) edge[-latex, bend left=10, black, opacity = 1.0] (5)
          (5) edge[-latex, bend left=10, black, opacity = 1.0] (4)
          (5) edge[-latex, bend left=10, black, opacity = 1.0] (3)
          (3) edge[-latex, bend left=10, red, opacity = 1.0] (5)
          (5) edge[-latex, black, opacity = 1.0] (1) ; % 5 5
\end{tikzpicture}
}
\caption{Illustration of directed ear decomposition: Starting with the 2-cycle in (a), we iteratively add ears, highlighted in red in each step, to obtain the  digraph in~(g).}\label{fig:eardec}
\end{figure}
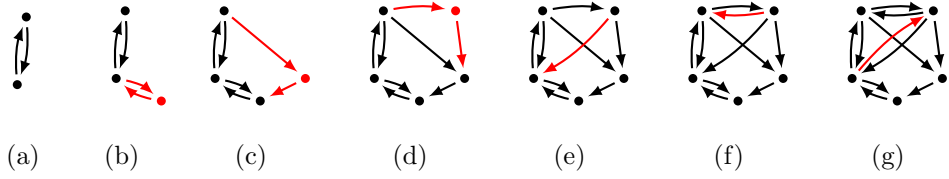

Let $k'$ be the number of cycles of $\vec S'$. We claim that $k'< k$, i.e., $\vec S$ contains more cycles than its subgraph $\vec S'$ does. To wit, if $u_{1} = u_{r}$, then $\vec P$ is a cycle of $\vec S$ but not of $\vec S'$. If $u_{1} \neq u_{r}$, then we let $\vec P'$ be a path in $\vec S'$ from $u_{r}$ to $u_{1}$. Concatenating $\vec P'$ with $\vec P$, we obtain a cycle in $\vec S$, which is not in $\vec S'$. This establishes the claim. 

Re-label the cycles of $\vec S$, if necessary, so that $\vec C_1,\ldots, \vec C_{k'}$ are the cycles of $\vec S'$, and $\vec C_{k'+1},\ldots, \vec C_{k}$ are the cycles of $\vec S$ but not of $\vec S'$. 
Each $\vec C_j$, for $j = k'+1,\ldots, k$, must contain the ear $\vec P$. 
Let 
$y^{(j)}:= c_j z_j$ for $j = 1,\ldots, k$ and   
$y':= \sum_{j = 1}^{k'} y^{(j)}$. 
Since all the $c_j$'s are positive, we have that 
$\supp(y^{(j)}) = V(\vec C_j)$. Since $\vec S'$ is strongly connected, every node of $\vec S'$ is contained in some cycle $\vec C_j$, for $j = 1,\ldots, k'$, and hence, $\supp(y') = V(\vec S')$. 
We then truncate $y'$ and the $y^{(j)}$'s by setting 
$$
\tilde y' := y' |_{\vec S'} \quad \mbox{and} \quad \tilde y^{(j)} := y^{(j)} |_{\vec C_j}, \quad \mbox{for all } j = k' + 1,\ldots, k.  
$$
For ease of notation, let
$$
\vec K:= \vec K_y(\vec S), \quad \vec K':= \vec K_{\tilde y'}(\vec S'), \quad \mbox{and} \quad \vec K^{(j)} := \vec K_{\tilde y^{(j)}}(\vec C_j)\quad \mbox{for all } j = k' + 1,\ldots, k.
$$
Since  $y = y' + \sum_{j = k' + 1}^{k} y^{(j)}$, one can embed simultaneously   
$\vec K'$ and $\vec K^{(j)}$, for $j = k' + 1,\ldots, k$, into $\vec K$. In other words, $\vec K$ contains these $(k - k' + 1)$ subgraphs whose node sets are pairwise disjoint.

Since $\vec S'$ is strongly connected and has~$k'$ cycles, for $k' < k$, and since $y'= \sum_{j = 1}^{k'} c_j z_j$ with the $c_j$'s positive, 
we can appeal to the induction hypothesis to obtain a Hamilton cycle $\vec H'$ of $\vec K'$. 
Through the embedding of $\vec K'$ into $\vec K$, we treat $\vec H'$ as a cycle of $\vec K$.  
Because $\vec S'$ contains the node $u_1$, there is a node $v'_{1}$ in $\vec H'$ such that $\pi(v'_1) = u_1$. We write $\vec H'$ explicitly as
\begin{equation}\label{eq:c'fromS'}
\vec H':= v'_1 \cdots v'_{n'} v'_1,   
\end{equation}
where $n' := \|y'\|_1 = |V(\vec K')|$. 

For each $j  = k'+1,\ldots, k$, we use the same arguments as in the base case to obtain a Hamilton cycle $\vec H^{(j)}$ of $\vec K^{(j)}$. Similarly, we treat $\vec H^{(j)}$ as a cycle of $\vec K$.  
Because $\vec H^{(j)}$ contains the ear $\vec P$ and, hence, the node $u_1$, there exists a node $v_{j,1}$ in $\vec H^{(j)}$ such that $\pi(v_{j,1}) = u_1$. We write $\vec H^{(j)}$ explicitly as
\begin{equation}\label{eq:cjfromdj}
\vec H^{(j)} := v_{j,1} \cdots v_{j,n_j} v_{j,1},
\end{equation}
where $n_j := \|y^{(j)}\| = |V(\vec K^{(j)})|$. 

Since $V(\vec K')$ and the $V(\vec K^{(j)})$ form a partition of $V(\vec K)$, their respective Hamilton cycles, namely, $\vec H'$ and the $\vec H^{(j)}$'s, form a cycle cover of $\vec K$. 
We will now use these cycles to construct a Hamilton cycle of $\vec K$. 
Since $\vec K$ is complete $\vec S$-partite and since the nodes $v'_1$ and $v_{j,1}$, for $j = k'+1,\ldots,k$, belong to $\pi^{-1}(u_1)$, we have that $v'_{n'}v_{k'+1,1}$, $v_{k,n_{k}} v'_1$, and $v_{j,n_j}v_{j+1,1}$, for $j = k'+1, \ldots, k-1$, are edges of $\vec K$. Thus,  
$$
\vec H := v'_1 \cdots v'_{n'}v_{k'+1,1} \cdots v_{k'+1, n_{k'+1}} v_{k'+2, 1} \cdots v_{k,n_{k}} v'_1
$$
is a desired Hamilton cycle of $\vec K$. 
\end{proof}

\section{Proof of Sufficiency of Conditions of $A$, $B$, and $C$}\label{sec:sufficiency}
In this section, we show that if a step-graphon $W$ satisfies Conditions $A$, $B$ (and $C$), then $W$ has the weak (strong) $H$-property.   

The condition that a graph has a cycle cover or a Hamilton cycle is monotone with respect to edge addition. Specifically, if $\vec G$ and $\vec G'$ are two graphs on the same node set, with $E(\vec G) \supseteq E(\vec G')$, then 
\begin{multline*}
G' \mbox{ has a cycle cover (Hamilton cycle)} \,\, \Rightarrow \,\,  G\mbox{ has a cycle cover (Hamilton cycle)}.
\end{multline*}
This monotonicity is carried over to graphons. Specifically if $W'$ and $W$ are two graphons, with $W' \leq W$, then
\begin{multline*}
\bP(\vec G_n\sim W' \mbox{ has a cycle cover (Hamilton cycle)}) \\
\leq \bP(\vec G_n\sim W \mbox{ has a cycle cover (Hamilton cycle)})
\end{multline*}
which implies that
$$
W' \mbox{ has weak (strong) } H\mbox{-property} \,\, \Rightarrow \,\,  
W \mbox{ has weak (strong) }H\mbox{-property}. 
$$
Thus, by Corollary~\ref{cor:loopfreegraphon}, we can assume that $W$ is loop free.

By Theorem~\ref{thm:constructHam}, for any $y\in \cX_0$, $\vec K_y$ has a cycle cover. If, further, $\vec S$ is strongly connected, then $\vec K_y$ has a Hamilton cycle. 
Denote by $\vec H_y$ the cycle cover (Hamilton cycle) as a subgraph of $\vec K_y$, and $\psi_y$ the embedding (i.e., a one-to-one graph homomorphism): 
\begin{equation}\label{eq:embeddingphi*}
\psi_y: \vec H_y \to \vec K_y.
\end{equation}   
We will show that if $W$ satisfies Conditions $A$ and $B$, then the same~$\vec H_y$ exists in $\vec G_n \sim W$ \aas. We make the statement precise below.

\begin{definition}
    Let $\vec G$ be an $\vec S$-partite graph, with $y(\vec G)\in \cX_0$. An embedding $\phi: \vec H_{y(\vec G)}\to \vec G$, if exists, is {\bf compatible with $\psi_{y(\vec G)}$} if $\pi\cdot \phi = \pi\cdot \psi_{y(\vec G)}$.
\end{definition}

Let $\vec{\mathcal{H}}$ be the set of $\vec S$-partite graphs $\vec G$ such that $y(\vec G)\in \cX_0$ and $\vec G$ admits an embedding $\phi: \vec{H}_{y(\vec G)} \to \vec{G}$, compatible with $\psi_{y(\vec G)}$. 
The  main result of the section is

\begin{theorem}\label{thm:almostsureembedding}
    Let $W$ be a loop-free step graphon. If $W$ satisfies Conditions $A$ and $B$, then 
    $\lim_{n\to\infty}\bP(\vec G_n \sim W \in \vec \cH) = 1$.
\end{theorem}

We take a two-step approach to establish the result: In Subsection~\ref{ssec:symmetrization}, we associate to $W$ a {\it symmetric} step-graphon $W^\rs$ and use it to sample {\it undirected} random graph $G_n\sim W^\rs$. The two random graphs $G_n$ and $\vec G_n$ relate to each other in the way that the probability of the event that $\vec G_n\in \vec \cH$ is bounded from below by the probability of the event that $H_{y(G_n)}$ is embeddable into $G_n$, where $H_y$ is the undirected counterpart of $\vec H_y$. This step allows for the use of the Blow-up Lemma, which we do in Subsection~\ref{ssec:blowuplemma}.

\subsection{Reduction by symmetrization}\label{ssec:symmetrization}

\begin{figure}[t]
\centering
\begin{subfigure}{.22\textwidth}
\centering
 \begin{tikzpicture}[scale=.45]
\filldraw [fill=Gray!30!black!80, draw=Gray!30!black!80] (0,0) rectangle (1/4,7/4); % 41
\filldraw [fill=Gray!70!black!20, draw=Gray!70!black!20] (1/4,7/4) rectangle (4/4,12/4); % 32
\filldraw [fill=Gray!70!black!70, draw=Gray!70!black!70] (4/4,12/4) rectangle (9/4,15/4); % 23
\filldraw [fill=Gray!30!black!30, draw=Gray!30!black!30] (9/4,0) rectangle (25/8,7/8); % 44
\filldraw [fill=Gray!30!black!30, draw=Gray!30!black!30] (25/8,7/8) rectangle (4,7/4); % 55
\filldraw [fill=Gray!10!black!20, draw=Gray!10!black!20] (9/4,7/4) rectangle (16/4,12/4); % 34
\filldraw [fill=Gray!60!black!70, draw=Gray!60!black!70] (1/4,15/4) rectangle (4/4,16/4); % 12
\filldraw [fill=Gray!30!black!10, draw=Gray!30!black!10] (4/4,0) rectangle (9/4,7/4); % 43

\draw [draw=black,very thick] (0,0) rectangle (4,4);

\node [above left] at (0,4) {$0$};
\node [left] at (0,0) {$1$};
\node [above] at (4,4) {$1$};
\end{tikzpicture}
\caption{$W$.}\label{sfig1:symgraphon}
\end{subfigure}
\begin{subfigure}{.22\textwidth}
\centering
\begin{tikzpicture}[scale=.45]
\filldraw [fill=Gray!30!black!80, draw=Gray!30!black!80] (0,0) rectangle (1/4,7/4); % 41
\filldraw [fill=Gray!30!black!80, draw=Gray!30!black!80] (9/4,15/4) rectangle (4,4); % 14
\filldraw [fill=Gray!70!black!10, draw=Gray!70!black!10] (1/4,7/4) rectangle (4/4,12/4); % 32
\filldraw [fill=Gray!70!black!10, draw=Gray!70!black!10] (4/4,12/4) rectangle (9/4,15/4); % 23
\filldraw [fill=Gray!70!black!5, draw=Gray!70!black!5] (9/4,7/4) rectangle (16/4,12/4); % 34
\filldraw [fill=Gray!70!black!5, draw=Gray!70!black!5] (4/4,0) rectangle (9/4,7/4); % 43
\filldraw [fill=Gray!60!black!70, draw=Gray!60!black!70] (1/4,15/4) rectangle (4/4,16/4); % 12
\filldraw [fill=Gray!60!black!70, draw=Gray!60!black!70] (0,12/4) rectangle (1/4,15/4); % 21
\filldraw [fill=Gray!30!black!5, draw=Gray!30!black!5] (9/4,0) rectangle (25/8,7/8); % 44
\filldraw [fill=Gray!30!black!5, draw=Gray!30!black!5] (25/8,7/8) rectangle (4,7/4); % 55

\draw [draw=black,very thick] (0,0) rectangle (4,4);

\node [above left] at (0,4) {$0$};
\node [left] at (0,0) {$1$};
\node [above] at (4,4) {$1$};
\end{tikzpicture}
\caption{$W^{\rs}$.}\label{sfig2:symgraphon}
\end{subfigure}
\begin{subfigure}{.23\textwidth}
\centering
    \begin{tikzpicture}[scale=0.77]
		\node [circle,fill=OliveGreen,inner sep=1.5pt,label=left:{$u_4$}] (4) at (-0.8, 0) {};
            \node [circle,fill=brown,inner sep=1.5pt,label=right:{$u_5$}] (5) at (0.8, 0) {};
		\node [circle,fill=red,inner sep=1.5pt,label=left:{$u_3$}] (3) at (-1.0, -1.5) {};
		\node [circle,fill=orange,inner sep=1.5pt,label=below:{$u_2$}] (2) at (0, -2) {};
		\node [circle,fill=cyan,inner sep=1.5pt,label=right:{$u_1$}] (1) at (1.0, -1.5) {};

  \path[draw,thick,shorten >=2pt,shorten <=2pt]
		 (1) edge (2) % 1 2
		 (2) edge (3) % 2 3
          (4) edge (1)
          (4) edge (3)
		 (4) edge (5)
          (5) edge (3)
          (5) edge (1); % 5 5
\end{tikzpicture}
\caption{$S$.}\label{sfig3:symgraphon}
\end{subfigure}
\begin{subfigure}{.22\textwidth}
\centering
    \begin{tikzpicture}[scale=0.72,rotate=-45]
    \tikzset{every loop/.style={}}
		\node [circle,fill=red,inner sep=1.2pt] (4) at (0.25, -0.25) {};
		\node [circle,fill=red,inner sep=1.2pt] (5) at (-0.25, 0.25) {}; 
            \node [circle,fill=red,inner sep=1.2pt] (3) at (-0.25, -0.25) {}; 
            %%%% 1
		\node [circle,fill=OliveGreen,inner sep=1.2pt] (6) at (-0.25, 1.5) {};
		\node [circle,fill=OliveGreen,inner sep=1.2pt] (7) at (0.25, 1.5) {};
		\node [circle,fill=brown,inner sep=1.2pt] (8) at (-0.25, 2) {};
		\node [circle,fill=brown,inner sep=1.2pt] (9) at (0.25, 2) {}; %%%% 2
		\node [circle,fill=cyan,inner sep=1.2pt] (14) at (2, 1.75) {};
		\node [circle,fill=orange,inner sep=1.2pt] (17) at (1.5, -0.25) {};
		\node [circle,fill=orange,inner sep=1.2pt] (18) at (2, 0) {};
	
  \path[draw, very thin, shorten >=2pt,shorten <=2pt]
		 (6) edge[very thin] (8)
		 (4) edge[very thin] (17)
		 (3) edge[very thin] (7)
          (5) edge[very thin] (9)
		 (5) edge[very thin] (18)
		 (14) edge[very thin] (18)
		 (9) edge[very thin] (14)
		 (7) edge[very thin] (14)
		 ;
\node at(0,-.51){};
\end{tikzpicture}
\caption{$G_{n} \sim W^\rs$.}\label{sfig5:samplegnfromws}
\end{subfigure}
\caption{Given the step-graphon $W$ in~(a), we follow~\eqref{eq:defhatw} to obtain the symmetric step-graphon $W^\rs$ in~(b). The undirected graph in~(c) is the skeleton graph of $W^\rs$ for the partition $\sigma = \frac{1}{16}(0,1,4,9,12.5,16)$. 
The undirected graph $G_{n}$ in~(d) is sampled from $W^\rs$, following steps $S1$ and $S'2$.  
}
    \label{fig:symgraphon}
\end{figure}
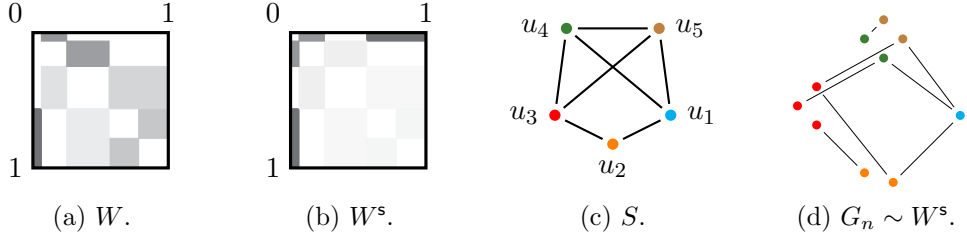

Let $\sigma = (\sigma_0,\ldots,\sigma_m)$ be a partition for $W$. 
Recall that $p_{ij}$ is the value of $W$ over $R_{ij} = [\sigma_{i-1}, \sigma_i)\times [\sigma_{j-1}, \sigma_j)$.    
We define a {\it symmetric} step-graphon $W^{\rs}$ as follows: 
For $1\leq i, j \leq m$ and for $(s, t)\in R_{ij}$, 
\begin{equation}\label{eq:defhatw}
W^{\rs}(s,t) := 
\begin{cases}
\max\{p_{ij}, p_{ji}\} & \mbox{if } p_{ij} p_{ji} = 0, \\
p_{ij} p_{ji} & \mbox{otherwise}. 
\end{cases} 
\end{equation}
We use $q_{ij}$ to denote the value of $W^{\rs}$ over $R_{ij}$ (and $R_{ji}$). 
See Figure~\ref{fig:symgraphon} for illustration.

To the step-graphon $W^{\rs}$ with partition $\sigma$, there corresponds the {\em undirected} graph $S$ on $m$ nodes, where a pair 
$(u_i,u_j)$ is an edge of $S$ if $q_{ij} > 0$.  
It follows from~\eqref{eq:defhatw} that $(u_i,u_j)$ is an edge of $S$ if and only if $\vec S$ contains either $u_iu_j$ or $u_ju_i$, or both. 
Since $\vec S$ does not have any self loop (as $W$ is loop free), neither does~$S$.

We use $W^{\rs}$ to sample an {\em undirected} graph $G_n$ on $n$ nodes as follows: First, follow step $S1$ (see Section~\ref{sec:introduction}) to obtain the coordinates $t_i$'s of the $n$ nodes. Then, follow
\begin{enumerate}    
\item[$\mathbf{S'2}$.] For each pair of two distinct nodes $v_i$ and $v_j$, place an {\em undirected} edge $(v_i,v_j)$ with probability $W^{\rs}(t_i, t_j)$. 
\end{enumerate}

It is clear that $G_n\sim W^{\rs}$ is $S$-partite (we omit the definition of $S$-partite graph as it is similar to the one for $\vec S$-partite graph). Denote by $K_y$ the complete $S$-partite graph, with $y(K_y) = y$.  
For any $y\in \cX_0$, let $H_y$ be the $S$-partite graph, as a subgraph of $K_y$, obtained from $\vec H_y$ by ignoring the orientations of its edges.  
If $\vec H_y$ is a cycle (with more than two nodes), then $H_y$ is an undirected cycle. 
If $\vec H_y$ is a cycle cover, then $H_y$ is a node-wise disjoint union of cycles and edges, where the edges correspond to the $2$-cycles in $\vec H_y$. 
With slight abuse of notation, we still use    
$\psi_y: H_y \to K_y$ to denote the embedding that sends the edges $(v_i, v_j)$ of $H_y$ to $(\psi_y(v_i), \psi_y(v_j))$.

Similarly, let $\cH$ be the set of all $S$-partite graphs $G$ such that $y(G) \in \cX_0$ and there exists an embedding $\phi: H_{y(G)}\to G$ compatible with $\psi_{y(G)}$, i.e., $\pi \cdot \phi = \pi\cdot \psi_{y(G)}$.   
The following result relates the event $\vec G_n\sim W\in \vec \cH$ to the event $G_n\sim W^\rs\in \cH$: 

\begin{proposition}\label{prop:undirectedcase}
    For any $n\in \N$, 
    $\bP(G_n \sim W^{\rs} \in \cH) \leq \bP(\vec G_n \in W \in \vec \cH)$.
\end{proposition}

We establish below the proposition. 
Given an $S$-partite graph $G$, we perform the following operation on its edge set to obtain an $\vec S$-partite digraph:
\begin{enumerate}

\item[$\mathbf{S'3}$.] For each edge $(v_i,v_j)$ of $G$, consider the following three cases: 

{\it Case 1: $\pi(v_i)\pi(v_j)\in \vec S$ and $\pi(v_j)\pi(v_i)\notin \vec S$.}  Replace $(v_i,v_j)$ with $v_iv_j$. 

{\it Case 2: $\pi(v_j)\pi(v_i)\in \vec S$ and $\pi(v_i)\pi(v_j)\notin \vec S$.} Replace $(v_i,v_j)$ with $v_jv_i$. 

{\it Case 3: $\pi(v_i)\pi(v_j), \pi(v_j)\pi(v_i)\in \vec S$.} Replace $(v_i,v_j)$ with $v_iv_j$ and $v_jv_i$. 
\end{enumerate}
We denote by $\vec G^{\rs}$ the resulting digraph.

Note that an embedding $\phi: H_y \to G$, with $y(G) = y$, does not necessarily induce an embedding $\phi: \vec H_y \to \vec G^{\rs}$; indeed, there may exist an edge $v_iv_j$ of $\vec H_y$ such that $\phi(v_i)\phi(v_j)$ is not an edge of $\vec G^\rs$. 
The following lemma shows that the induced embedding always exists if $\phi$ is compatible with $\psi_y$:

\begin{lemma}\label{lem:reverseinduce}
Let $G \in \cH$ and $\phi: H_{y(G)}\to G$ be an embedding compatible with $\psi_{y(G)}$. Then, $\phi$ induces an embedding of $\vec H_{y(G)}$ to $\vec G^{\rs}$.  In particular, 
$G \in \cH$ if and only if $\vec G^\rs \in \vec \cH$.
\end{lemma}

\begin{proof}
Within the proof, we will simply write $\psi$ by omitting its sub-index. 
We show that if $v_iv_j\in E(\vec H_{y(G)})$, then $\phi(v_i)\phi(v_j)\in E(\vec G^\rs)$. 
Since $\phi$ is compatible with $\psi$, we have that
$u_i := \pi\cdot \psi(v_i) = \pi  \cdot \phi(v_i)$ and $u_j := \pi\cdot \psi(v_j) = \pi\cdot \phi(v_j)$. 
Since $\pi\cdot \psi: \vec H_y\to \vec S$ is a graph homomorphism, $u_iu_j$ is an edge of $\vec S$. Also, since $\phi: H_y\to G$ is a graph homomorphism, $(\phi(v_i), \phi(v_j))$ is an edge of $G$.  
Then, by the operation given in the step $S'3$, we conclude that $\phi(v_i)\phi(v_j)$ is an edge of $\vec G^{\rs}$. 
\end{proof}

With slight abuse of notation, we denote by $\vec G^{\rs}_n\sim W^{\rs}$  the random digraph on $n$ nodes obtained by following the steps $S1$, $S'2$, and $S'3$. 
An immediate consequence of Lemma~\ref{lem:reverseinduce} is then the following: 

\begin{lemma}\label{lem:intermediatestep0}
For any $n\in \N$, 
\begin{equation}\label{eq:intermediatestep0}
\bP(G_n \sim W^\rs \in \cH ) = \bP(\vec G^\rs_n \sim W^\rs \in \vec \cH).
\end{equation}
\end{lemma}

The following lemma relates the event $\vec G^\rs_n \sim W^\rs\in \cH$ to the event $\vec G_n\sim W\in \cH$, and completes the proof of  Proposition~\ref{prop:undirectedcase}. 

\begin{lemma}\label{lem:intermediatestep}
For any $n\in \N$,  
$$
\bP(\vec G^{\rs}_n \sim W^{\rs} \in \vec \cH) \leq \bP(\vec G_n \sim W \in \vec \cH). 
$$
\end{lemma}

\begin{proof}
Given an arbitrary $\vec S$-partite graph $\vec G_n$, we let $\vec G^*_n$ be obtained by removing certain edges out of $\vec G_n$ as specified below:
\begin{enumerate}
\item[$\mathbf{S3}$.] We remove an edge $v_iv_j$ of $\vec G_n$ if both of the following two conditions hold:  

{\it 1:} Both $\pi(v_i)\pi(v_j)$ and $\pi(v_j)\pi(v_i)$ are edges of $\vec S$.   

{\it 2:}  $v_jv_i$ is not an edge of $\vec G_n$.
\end{enumerate}

Denote by $\vec G_n^*\sim W^*$ the random digraph obtained by following the steps $S1$, $S2$, and $S3$. 
Let $u_i, u_j\in V(\vec S)$ be such that $u_iu_j, u_ju_i\in E(\vec S)$.  
It is clear that for two distinct nodes $v_i\in \pi^{-1}(u_i)$ and $v_j\in \pi^{-1}(u_j)$, the probability that $\vec G_n\sim W$ has both edges $v_iv_j$ and $v_j v_i$ is $p_{ij}p_{ji} = q_{ij}$, so 
$$
\bP(v_iv_j\in \vec G^*_n \mbox{ and } v_jv_i\in \vec G^*_n) = q_{ij} \quad \mbox{and} \quad \bP(v_iv_j\notin \vec G^*_n \mbox{ and } v_jv_i\notin \vec G^*_n) = 1 - q_{ij}. 
$$
It follows that the two sampling procedures, namely, the one ($S1$-$S'2$-$S'3$) for sampling $\vec G^{\rs}_n \sim W^{\rs}$ and the other ($S1$-$S2$-$S3$) for sampling $\vec G_n^*\sim W^*$, are equivalent to each other. It follows that
\begin{equation}\label{eq:intermediatestep1}
    \bP(\vec G^{\rs}_n \sim W^{\rs} \in \vec \cH) = \bP(\vec G_n^* \sim W^* \in \vec \cH). 
\end{equation}

The condition that an $\vec S$-partite graph belongs to $\vec \cH$ is monotone with respect to edge addition. 
Since $\vec G^*_n$ is obtained from $\vec G_n$ by removing edges, $\vec G^*_n\in \vec \cH$ implies $\vec G_n \in \vec \cH$. Thus,  
\begin{equation}\label{eq:intermediatestep2}
\bP(\vec G_n^* \sim W^* \in \vec \cH) \leq \bP(\vec G_n \sim W \in \vec \cH). 
\end{equation}
The lemma then follows from~\eqref{eq:intermediatestep1} and~\eqref{eq:intermediatestep2}. 
\end{proof}

\subsection{On the use of the Blow-up Lemma}\label{ssec:blowuplemma}

Let $\cU$ be the open neighborhood of $x^*$ in $\cX$, introduced at the beginning of Section~\ref{sec:constructHam}. 
Since $x^*\in \rint \cX$ and since $x(G_n)$ converges to $x^*$ \aas, it holds that   
$y(G_n) \in \cX_0$ \aas. 
We show below that 
\begin{equation}\label{eq:alsmotsureembddingH}
\lim_{n\to\infty} \bP(G_n \sim W^{\rs} \in \cH \mid y(G_n) \in \cX_0) = 1. 
\end{equation}  
In words, we show that if $y(G_n)\in \cX_0$, then \aas~there exists an embedding $\phi: H_{y(G_n)} \to G_n$, compatible with $\psi_{y(G_n)}: H_{y(G_n)}\to K_{y(G_n)}$. Note that if~\eqref{eq:alsmotsureembddingH} holds, then by Proposition~\ref{prop:undirectedcase}, we have that
$\lim_{n\to\infty}\bP(\vec G_n \sim W \in\vec H) = 1$, i.e., Theorem~\ref{thm:almostsureembedding} holds, which will then complete the proof of~\eqref{eq:posevent} and~\eqref{eq:posevent1}.  

The proof of~\eqref{eq:alsmotsureembddingH} relies on the use of the {Blow-up Lemma}, which we recall below.  
Let $G$ be an arbitrary undirected graph. For two disjoint subsets $X$ and $Y$ of $V(G)$, let $e(X,Y)$ be the number of edges between $X$ and $Y$. We need the following definition:

\begin{definition}[Super-regular pair]
    Let $G$ be an undirected graph, and $A$, $B$ be two disjoint subsets of $V(G)$. The pair $(A,B)$ is {\bf $(\epsilon,\delta)$-super-regular} if 
    \begin{equation}\label{eq:edgedensitycondition}
    e(X,Y) > \delta |X||Y|, \mbox{ for any } X\subseteq A \mbox{ and } Y\subseteq B, \mbox{ with } |X| > \epsilon |A| \mbox{ and } |Y| > \epsilon |B|,  
    \end{equation}
    and, moreover,   
    \begin{equation}\label{eq:degreecondition}
    e(a, B) > \delta |B| \quad \mbox{for any } a\in A, \quad \mbox{and} \quad e(b, A) > \delta |A| \quad \mbox{for any } b\in B.
    \end{equation}
\end{definition}

We extend the above definition to the $S$-partite graphs:

\begin{definition}[Super-regular $S$-partite graphs] 
Let $S$ be an undirected graph, without self-loops, on $m$ nodes.  
An $S$-partite graph $G$, with $y(G)\in \N^m$, is {\bf $(\epsilon,\delta)$-super-regular} if for any two distinct nodes $u_i, u_j\in V(S)$,  $(\pi^{-1}(u_i), \pi^{-1}(u_j))$ is $(\epsilon,\delta)$-super-regular.
\end{definition}

For an arbitrary graph $H$, let $\Delta(H)$ be the degree of $H$ (i.e., the maximum of the degrees of its nodes).  
We reproduce below the Blow-up Lemma~\cite{komlos1997blow}: 

\begin{lemma}[Blow-up Lemma]\label{lem:blowuplemma} 
Let $S$ be an undirected graph, without self-loops, on $m$ nodes. Then, given parameters $\delta > 0$ and $\Delta\in \N$, there exists an $\epsilon  = \epsilon(\delta,\Delta, m) > 0$ such that for any $y\in \N^m$, the following holds: If $H$ is an undirected graph with $\Delta(H) \leq \Delta$ and if there is an embedding $\psi: H \to K_y(S)$, then for any $(\epsilon, \delta)$-super-regular $S$-partite graph $G$, with $y(G) = y$, there is an embedding $\phi: H \to G$, compatible with $\psi$.    
\end{lemma}

We now return to the proof of~\eqref{eq:alsmotsureembddingH}. 
For any $y\in \cX_0$, we let $H_y$ be given as in the previous subsection. 
As argued earlier, $H_y$ is either a cycle or a node-wise disjoint union of cycles and possibly edges.  
Thus, 
$\Delta(H_y) \leq 2$, for all $y\in \cX_0$.  
Also, by Theorem~\ref{thm:constructHam}, there is an embedding $\psi_y: H_y\to K_{y}$ for all $y\in \cX_0$.  
Let 
\begin{equation}\label{eq:defdelta}
\delta := \frac{1}{2} \min \{q_{ij} \mid (u_i,u_j) \in E(S)\}.
\end{equation}
and $\epsilon := \epsilon(\delta, 2, m) > 0$ be given as in the statement of Lemma~\ref{lem:blowuplemma}. 
It remains to show that \aas~$G_n\sim W^\rs$ is $(\epsilon,\delta)$-super-regular. 

\begin{proposition}\label{prop:superregular}
For any $\epsilon > 0$, 
$$
\lim_{n\to\infty} \bP(G_n \sim W^{\rs} \mbox{ is } (\epsilon,\delta)\mbox{-super-regular}) = 1. 
$$
\end{proposition}

The proof of the proposition uses standard arguments in random graph theory. 
For completeness of presentation, we include it in Appendix~\ref{sec:blowuplemma}.  \hfill{\qed}

\printbibliography

\appendix 

\section{Proof of Proposition~\ref{prop:equivalenceofpartition}}\label{sec:partitionsame}

We say that $\sigma'$ is a {\it refinement} of $\sigma$ if $\sigma'$ contains $\sigma$ as a subsequence. Furthermore, $\sigma'$ is a {\it one-step refinement} of $\sigma$ if $\sigma'$ contains one more element than $\sigma$ does. It is clear that any refinement can be obtained by iterating one-step refinements.    
Note that for any two arbitrary partitions $\sigma$ and $\sigma'$, there exists a partition $\sigma''$ as a refinement of both $\sigma$ and $\sigma'$. 
The arguments above then imply that to establish Proposition~\ref{prop:equivalenceofpartition}, it suffices to prove for the case where $\sigma'$ is a one-step refinement of $\sigma$.  

\begin{figure}[h]
\centering
\begin{subfigure}{.45\textwidth}
\centering
    \begin{tikzpicture}[scale=1]
    \tikzset{every loop/.style={}}
		\node [circle,fill=black,inner sep=1.2pt,label=left:{$u_4$}] (4) at (-0.8, 0) {};
            \node [circle,fill=black,inner sep=1.2pt,label=right:{$u_5$}] (5) at (0.8, 0) {};
		\node [circle,fill=black,inner sep=1.2pt,label=left:{$u_3$}] (3) at (-1.0, -1.5) {};
		\node [circle,fill=black,inner sep=1.2pt,label=below:{$u_2$}] (2) at (0, -2) {};
		\node [circle,fill=black,inner sep=1.2pt,label=right:{$u_1$}] (1) at (1.0, -1.5) {};

  \path[draw,thick,shorten >=2pt,shorten <=2pt, -latex]
		 (1) edge[-latex, black, opacity = 1] (2) % 1 2
		 (2) edge[-latex, bend left=15, black, opacity = 1] (3) % 2 3
          (3) edge[-latex, bend left=15, black, opacity = 1] (2)
          (4) edge[-latex, black, opacity = 1] (1)
          (4) edge[-latex, bend left=12, black, opacity = 1] (3)
          (3) edge[-latex, bend left=12, black, opacity = 1] (4)
		 (4) edge[-latex, bend left=10, black, opacity = 1.0] (5)
          (5) edge[-latex, bend left=10, black, opacity = 1.0] (4)
          (5) edge[-latex, bend left=10, black, opacity = 1.0] (3)
          (3) edge[-latex, bend left=10, black, opacity = 1.0] (5)
          (5) edge[-latex, black, opacity = 1.0] (1) % 5 5
          (4) edge[loop above, opacity = 1] (4)
          (5) edge[loop above, black] (5);
\end{tikzpicture}
\caption{Digraph $\vec S'$.}\label{sfig2:refiness'}
\end{subfigure}
 \begin{subfigure}{.45\textwidth}
\centering
\begin{tikzpicture}[scale=1]

		\node [circle,fill=black,inner sep=1pt,label=above:{$u'_1$}] (1) at (0, 1.2) {};
		\node [circle,fill=black,inner sep=1pt,label=above:{$u'_2$}] (2) at (1, 1.2) {};
		\node [circle,fill=black,inner sep=1pt,label=above:{$u'_3$}] (3) at (2, 1.2) {};
		\node [circle,fill=black,inner sep=1pt,label=above:{$u'_4$}] (4) at (3, 1.2) {};
        \node [circle,fill=black,inner sep=1pt,label=above:{$u'_5$}] (5) at (4, 1.2) {};
		\node [circle,fill=black,inner sep=1pt,label=below:{$u''_1$}] (11) at (0, -0.5) {};
		\node [circle,fill=black,inner sep=1pt,label=below:{$u''_2$}] (12) at (1, -0.5) {};
		\node [circle,fill=black,inner sep=1pt,label=below:{$u''_3$}] (13) at (2, -0.5) {};
		\node [circle,fill=black,inner sep=1pt,label=below:{$u''_4$}] (14) at (3, -0.5) {};
        \node [circle,fill=black,inner sep=1pt,label=below:{$u''_5$}] (15) at (4, -0.5) {};
				 
		 \path[draw,thick,shorten >=2pt,shorten <=2pt]
		 (4) edge [black, opacity = 1.0] (11)
		 (1) edge [black, opacity = 1.0] (12)
		 (2) edge [black, opacity = 1] (13)
		 (3) edge [black, opacity = 1] (12)
          (3) edge [black, opacity = 1] (14)
          (4) edge [black, opacity = 1] (13)
          (4) edge [black, opacity = 1] (14)
          (4) edge [black] (15)
          (5) edge [black] (14)
          (3) edge [black] (15)
          (5) edge [black] (13)
		 (5) edge [black] (11)
          (5) edge [black] (15);
\end{tikzpicture}
\caption{Bipartite graph $B_{\vec S'}$.}\label{sfig2:refinebb'}
\end{subfigure}   
\caption{The digraph $\vec S'$ is the one-step refinement of the graph $\vec S$ shown in Figure~\ref{sfig:bipar1} on node $u_4$. The the bipartite graph  $B_{\vec S'}$ is associated with $\vec S'$, respectively.
 }\label{fig:onesteprefinement} 
\end{figure}
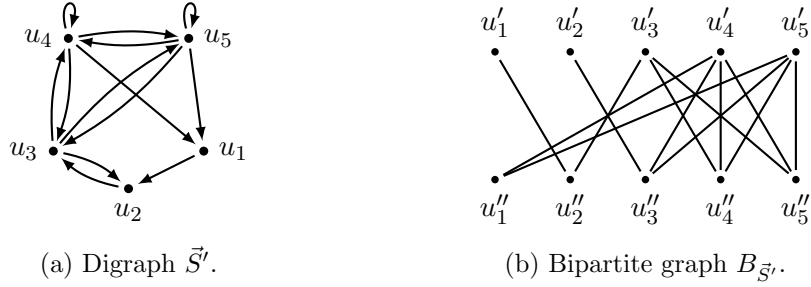

Let $\sigma = (\sigma_0,\ldots, \sigma_{m-1}, \sigma_*)$, with $\sigma_0 = 0$ and $\sigma_* = 1$. We assume, without loss of generality, that $\sigma'$ is obtained from $\sigma$  by inserting an element $\sigma_m$ between $\sigma_{m-1}$ and~$\sigma_*$, i.e.,  
$\sigma' = (\sigma_0,\ldots, \sigma_{m-1},\sigma_m, \sigma_*)$. 
Then, the following hold for $x'^*$, $\vec S'$, and $B_{\vec S'}$ (also, see Figure~\ref{fig:onesteprefinement} for illustration):
\begin{enumerate}
\item Let $\hat x^* := (x^*; 0) \in \R^{m+1}$. Then, 
\begin{equation}\label{eq:defbeta}
x'^* = \hat x^* +  (1 - \sigma_m) (e_{m + 1}  - e_m). 
\end{equation}
\item The skeleton graph $\vec S'$ can be obtained from $\vec S$ by adding the new node $u_{m+1}$ and the new edges incident to $u_{m+1}$: 
\begin{multline}\label{eq:edgesetfors'refine}
E(\vec S')  = E(\vec S) \cup \{u_i u_{m+1} \mid u_iu_m\in E(\vec S)\} \\ \cup \{u_{m+1} u_j \mid u_mu_j \in E(\vec S)\}   
\cup \{u_{m+1}u_{m+1} \mbox{ if } u_mu_m\in E(\vec S) \}. 
\end{multline}
We call $\vec S'$ the {\em one-step refinement} of $\vec S$ on node~$u_m$. 

\item Correspondingly, the node and the edge sets of $B_{\vec S'}$ are given by
$$
V'(B_{\vec S'})  = V'(B_{\vec S}) \cup \{u'_{m+1}\}, \quad  
V''(B_{\vec S'})  = V''(B_{\vec S}) \cup \{u''_{m+1}\}, 
$$
and
\begin{equation}\label{eq:edgesetforbs'refine}
\begin{aligned}
E(B_{\vec S'})  = E(B_{\vec S}) & \cup \{(u'_i, u''_{m+1}) \mid (u'_i,u''_m)\in E(B_{\vec S})\} \\
& \cup \{(u'_{m+1}, u''_j) \mid (u'_m, u''_j) \in E(B_{\vec S})\} \\
& \cup \{(u'_{m+1}, u''_{m+1}) \mbox{ if } (u'_m, u''_m) \in E(B_{\vec S}) \}.
\end{aligned}
\end{equation}
\end{enumerate}
We now establish the three items of Proposition~\ref{prop:equivalenceofpartition}:

\subsection{Proof of item 1}
We define the graph homomorphism $\theta: \vec S' \to \vec S$ as
$$
\theta(u_i):= 
\begin{cases}
u_i & \mbox{if } 1\leq i \leq m, \\
u_{m} & \mbox{if } i = m + 1.
\end{cases}
$$
It follows from~\eqref{eq:edgesetfors'refine} that 
$u_i u_j \in E(\vec S')$ if and only if  $\theta(u_i)\theta(u_j) \in E(\vec S)$.  
Thus, 
\begin{equation}\label{eq:walksofss'}
\vec P'=u_{i_1}\cdots u_{i_\ell} \mbox{ is a walk of } \vec S' \,\, \Leftrightarrow \,\,  \theta(\vec P'):=\theta(u_{i_1})\cdots \theta(u_{i_\ell}) \mbox{ is a walk of } \vec S.  
\end{equation}

\xc{Proof that $\vec S'$ is strongly connected $\Rightarrow$ $\vec S$ is strongly connected.} Let $u_i$ and $u_j$ be two distinct nodes of $\vec S$. We pick nodes $u_{i'}\in \theta^{-1}(u_i)$ and $u_{j'}\in \theta^{-1}(u_j)$. Since $\vec S'$ is strongly connected, there is a path $\vec P'$ of $\vec S'$ from $u_{i'}$ to $u_{j'}$. By~\eqref{eq:walksofss'}, we have that $\theta(\vec P')$ is a walk of $\vec S$ from $u_i$ to $u_j$. 

\xc{Proof that $\vec S$ is strongly connected $\Rightarrow$ $\vec S'$ is strongly connected.} Let $u_{i'}$ and $u_{j'}$ be two distinct nodes of $\vec S'$. We first consider the case where $u_i := \theta(u_{i'})$ and $u_j := \theta(u_{j'})$ are two distinct nodes.  
In this case, there is a path $\vec P = u_{i_1}\cdots u_{i_\ell}$, with $u_{i_1} = u_i$ and $u_{i_\ell} = u_j$, of $\vec S$ from $u_i$ to $u_j$. We pick nodes $u_{i'_j}\in \theta^{-1}(u_{i_j})$, for $j = 2,\ldots, \ell-1$. Then, by~\eqref{eq:walksofss'}, $u_{i'} u_{i'_2} \cdots u_{i'_{\ell - 1}}u_{j'}$ is a path of $\vec S'$ from $u_{i'}$ to $u_{j'}$.   
We now assume that $u_i = u_j$. Since $\vec S$ has at least $2$ nodes, there exists a node $u_k$ of $\vec S$ such that $u_k \neq u_i$. Let $\vec P_{1} = u_{i_1} \cdots u_{i_{\ell_1}}$  (resp., $\vec P_{2} = u_{i_{\ell_1}} u_{i_{\ell_1 + 1}}\cdots u_{i_{\ell}}$) be a path of $\vec S$ from $u_i$ to $u_k$ (resp., from $u_k$ to $u_i$), where $u_{i_1} = u_{i_\ell} = u_i$ and $u_{i_{\ell_1}} = u_k$. Concatenating $\vec P_1$ and $\vec P_2$, we obtain a closed walk. Pick nodes $u_{i'_j}\in \theta^{-1}(u_{i_j})$, for $j = 2,\ldots, \ell-1$. Using again~\eqref{eq:walksofss'}, we conclude that $u_{i'}u_{i'_2} \cdots u_{i'_{\ell - 1}} u_{j'}$ is a walk of $\vec S'$ from $u_{i'}$ to $u_{j'}$.  \hfill{\qed}

\subsection{Proof of item~2} 
Let $\vec S_1,\ldots, \vec S_q$ be the SCCs of $\vec S$, and $\vec S'_1,\ldots, \vec S'_{q'}$ be the SCCs of $\vec S'$. Without loss of generality, we assume that $u_m\in V(\vec S_q)$. 
By Lemma~\ref{lem:connectednessofB}, it suffices to show that 
\begin{equation}\label{eq:reducetoconnectBitem2}
B_{\vec S_1},\ldots, B_{\vec S_q} \mbox{ are connected} \quad  \Longleftrightarrow \quad 
B_{\vec S'_1},\ldots, B_{\vec S'_{q'}} \mbox{ are connected}.
\end{equation}

If $\vec S_q$  comprises  the single node $u_m$ without self-loop, then $\vec S'$ has $(q + 1)$ SCCs $\vec S'_1,\ldots, \vec S'_{q+1}$, where  
$\vec S'_p := \vec S_p$, for $p = 1,\ldots, q$, and $\vec S'_{q+1}$ comprises the single node $u_{m+1}$ without self-loop. 
But then, the bipartite graph $B_{\vec S_q}$ has two nodes $u'_m$ and $u''_m$, without the edge $(u'_m,u''_m)$, so $B_{\vec S_q}$ is disconnected. The same applies to $B_{\vec S'_q}$ and $B_{\vec S'_{q+1}}$. Thus, for one of the two sides of~\eqref{eq:reducetoconnectBitem2} to hold, we must have that either $\vec S_q$ has at least two nodes, or, $\vec S_q$ comprises the single node $u_m$ with self-loop.  
It follows that $\vec S'$ has $q$ SCCs 
$\vec S'_1,\ldots, \vec S'_q$, where $\vec S'_p := \vec S_p$ for all $p = 1,\ldots, q-1$, and $\vec S'_q$ is a one-step refinement of $\vec S_q$ on $u_m$. 

For convenience but without loss of generality, we now assume that $\vec S$ is itself strongly connected (and has a self-loop if it only has a single node). 
We show below that $B_{\vec S}$ is connected if and only if $B_{\vec S'}$ is.
The arguments are similar to those for proving item~1 of Proposition~\ref{prop:equivalenceofpartition}. With slight abuse of notation, we now let $\theta: B_{\vec S'} \to B_{\vec S}$ be the graph homomorphism defined as 
$$
\theta(u'_i):= 
\begin{cases}
u'_i & \mbox{if } 1\leq i \leq m \\
u'_{m} & \mbox{if } i = m + 1,
\end{cases} \quad 
\mbox{and} \quad
\theta(u''_i):= 
\begin{cases}
u''_i & \mbox{if } 1\leq i \leq m \\
u''_{m} & \mbox{if } i = m + 1.
\end{cases}
$$
It follows from~\eqref{eq:edgesetforbs'refine} that 
$(u'_i,u''_j) \in E(B_{\vec S'})$ if and only if  $(\theta(u'_i),\theta(u''_j)) \in E(B_{\vec S})$. 
Thus, $u'_{i_1}u''_{j_1} \cdots u'_{i_\ell} u''_{j_\ell}$ is a walk of  $B_{\vec S'}$ if and only if 
$\theta(u'_{i_1})\theta(u''_{j_1}) \cdots \theta(u'_{i_\ell}) \theta(u''_{j_\ell})$ is a walk of $B_{\vec S}$. 
This completes the proof.  \hfill{\qed}

\subsection{Proof of item~3}
We consider two cases: (1) $u_m$ does not have a self-loop and (2) $u_m$ has a self-loop.   
\subsubsection{Case 1: $u_m$ does not have a self-loop}
Let $\vec C_1, \ldots, \vec C_\ell$, for $\ell \leq k$, be the cycles of $\vec S$ that contain $u_m$. Every such cycle $\vec C_p$, for $1 \leq p \leq \ell$, induces  two different cycles of $\vec S'$: One is $\vec C_{p, 1}:= \vec C_{p}$ and the other  is obtained by replacing the node $u_m$ in $\vec C_p$ with $u_{m+1}$, which we denote by $C_{p, 2}$. The set of cycles of $\vec S'$ is thus given by 
$$\{\vec C_{p, i} \mid 1 \leq p \leq \ell \mbox{ and } 1\leq i \leq 2\} \cup \{\vec C_q \mid \ell + 1 \leq q\leq k\}.$$
Let $z'_{p, i}$ and $z'_q$ be the node-cycle incidence vectors of $\vec S'$ corresponding to $\vec C_{p, i}$ and $\vec C_q$, respectively. Then,  
\begin{equation}\label{eq:defz'forcase1}
z'_{p, 1} = \hat z_p, \quad z'_{p, 2} = \hat z_p - e_m + e_{m+1}, \quad \mbox{and } z'_q = \hat z_q,
\end{equation}
where we recall that $\hat z_j= (z_j; 0)$.  

\xc{Proof that $x^*\in \cX \Rightarrow x'^*\in \cX'$ {\em ($x^*\in \rint \cX \Rightarrow x'^*\in \rint \cX'$)}.}
We write $x^* = \sum_{j= 1}^k c_j z_j$ with $c_j \geq 0$. Since $\vec C_1,\ldots, \vec C_\ell$ are the cycles that contain $u_m$, we have that
\begin{equation}\label{eq:sumcp}
\sum_{p = 1}^\ell c_p = x^*_{m} = (1 - \sigma_{m-1}).
\end{equation}
Now, let
\begin{equation}\label{eq:defc'forcase1}
c'_{p, 1} := \frac{\sigma_m - \sigma_{m-1}}{1 - \sigma_{m-1}} c_p, \quad c'_{p, 2} := \frac{1 - \sigma_m}{1 - \sigma_{m-1}} c_p, \quad \mbox{and } 
c'_q := c_q.  
\end{equation}
Note  that $c'_{p, 1} + c'_{p, 2} = c_p$ for all $1 \leq p \leq \ell$. 
Then, 
\begin{align}
x'^* & = \hat x^* + (1 - \sigma_m) (e_{m + 1} - e_m) \notag\\
& = \sum_{j = 1}^k c_j \hat z_j + (1 - \sigma_m) (e_{m + 1} - e_m) \notag\\
& = \sum_{p = 1}^\ell \left [ \sum_{i = 1}^2 c'_{p, i} \right ]\hat z_{p} + \sum_{q = \ell + 1}^k c'_q \hat z_q + (1 - \sigma_m)(e_{m + 1} - e_m) \notag \\
& = \sum_{p = 1}^\ell \sum_{i = 1}^2 c'_{p, i} z'_{p, i} + \sum_{q = \ell + 1}^k c'_q z'_q + \left [ (1 - \sigma_{m}) - \sum_{p = 1}^\ell c'_{p, 2} \right ] (e_{m + 1} - e_m) \notag \\
& =  \sum_{p = 1}^\ell \sum_{i = 1}^2 c'_{p, i} z'_{p, i} + \sum_{q = \ell + 1}^k c'_q z'_q + (1 - \sigma_{m}) \left [ 1 - \frac{1}{ 1- \sigma_{m-1}} \sum_{p = 1}^\ell c_p \right ] (e_{m + 1} - e_m) \notag \\
& = \sum_{p = 1}^\ell \sum_{i = 1}^2 c'_{p, i} z'_{p, i} + \sum_{q = \ell + 1}^k c'_q z'_q, \label{eq:positivecombforx'case1}
\end{align}
where the first equality follows from~\eqref{eq:defbeta}, the fourth equality follows from~\eqref{eq:defz'forcase1}, the fifth equality follows from~\eqref{eq:defc'forcase1}, and the last equality follows from~\eqref{eq:sumcp}.  
By~\eqref{eq:positivecombforx'case1}, $x'^*\in \cX'$. 
If, further, the coefficients $c_j$'s are positive (which holds if $x^*\in \rint \cX$), then by~\eqref{eq:defc'forcase1} the $c'_{p, i}$'s and the $c'_q$'s are positive as well and hence, $x'^*\in \rint \cX'$.  

\xc{Proof that $x'^*\in \cX' \Rightarrow x^*\in \cX$ {\em ($x'^*\in \rint \cX' \Rightarrow x^*\in \rint \cX$)}.} We write 
$$x'^* = \sum_{p = 1}^\ell \sum_{i = 1}^2 c'_{p, i} z'_{p, i} + \sum_{q = \ell + 1}^k c'_q z'_q,$$
where the $c'_{p, i}$'s and the $c_q$'s are nonnegative. 
For $p = 1,\ldots, \ell$ and for $q = \ell + 1, \ldots, k$, we define 
\begin{equation}\label{eq:defcjforcase1}
c_p:= c'_{p, 1} + c'_{p, 2} \quad \mbox{and} \quad c_q:= c'_q. 
\end{equation}
Let $J\in \R^{m\times (m + 1)}$ be defined as follows: 
$$
J := 
\begin{bmatrix}
1 & & & \\
& \ddots & & \\
&  & 1 & 1
\end{bmatrix}.
$$
It follows from~\eqref{eq:defz'forcase1} and~\eqref{eq:defbeta} that  
$z_p = J z'_{p, i}$, $z_q = J z'_q$, and $x^* = J x'^*$. 
Thus, 
\begin{multline*}
x^*  = J x'^*  = \sum_{p = 1}^\ell \sum_{i = 1}^2 c'_{p, i} J z'_{p, i} + \sum_{q = \ell + 1}^k c'_q J z'_q  = \sum_{p = 1}^\ell \left [ \sum_{i = 1}^2 c'_{p, i} \right ] z_{p} + \sum_{q = \ell + 1}^k c'_q z_q = \sum_{j = 1}^k c_j z_j, 
\end{multline*}
which shows that $x^*\in \cX$. By~\eqref{eq:defcjforcase1}, if the $c'_{p, i}$'s and the $c'_q$'s are positive, then so are the $c_j$'s, which implies that if $x'^*\in \rint \cX'$, then $x^*\in \rint \cX$. \hfill{\qed}

\subsubsection{Case 2: $u_m$ has a self-loop}
We again let $\vec C_1, \ldots, \vec C_\ell$, for $\ell \leq k$, be the cycles of $\vec S$ that contain $u_m$, with $\vec C_1 = u_m u_m$ the self-loop. 
The self-loop $\vec C_1$ induces three cycles of $\vec S'$, which are  
$\vec C_{1, 1} = u_mu_m$, $\vec C_{1, 2} = u_{m + 1} u_{m + 1}$, and $\vec C_{1, 3} = u_m u_{m + 1} u_m$. 
As argued at the beginning of Subsection~\ref{ssec:proofthatx'inY'noloop}, 
each cycle $\vec C_p$, for $2\leq p \leq \ell$, induces four different cycles: 
$\vec C_{p, 1} := \vec C_p$ and $\vec C_{p, 2}$, $\vec C_{p, 3}$, $\vec C_{p, 4}$ are obtained by replacing $u_m \in \vec C_p$ with $u_{m+1}$, $u_mu_{m+1}$, $u_{m+1}u_m$, respectively. 

Let $z'_{1,i}$, for $i = 1,\ldots, 3$, be the node-cycle incidence vectors corresponding to $\vec C_{1,i}$, which are given by 
$z'_{1,1} = e_m$, $z'_{1,2} = e_{m+1}$, and $z'_{1,3} = e_m + e_{m + 1}$. 
Let $z'_{p,i}$, for $2\leq p \leq \ell$ and $1\leq i \leq 4$, be the node-cycle incidence vectors corresponding to $\vec C_{p,i}$, as given in~\eqref{eq:relatez'toz}. Note, in particular, that
$$
\begin{aligned}
z'_{1,3} & = z'_{1,1} + z'_{1,2}, \\
z'_{p,3} & = z'_{p,4}  = z'_{p,1} + z'_{1,2} = z'_{p,2} + z'_{1,1}, \quad \mbox{for } p = 2,\ldots, \ell,\\
\end{aligned}
$$
which implies that none of the vectors $z'_{1,3}$, $z'_{p,3}$, and $z'_{p,4}$ is an extremal generator of $\cX'$, and can thus be suppressed in the nonnegative (positive) combination of $x'^*$. 
The same arguments in the previous case can be used to establish the current case.    \hfill{\qed}

\section{On graphons with symmetric support}\label{sec:symmetricgraphon}

Let $W$ be a step-graphon with symmetric support, i.e., $W(s, t)\neq 0$ if and only if $W(t,s) \neq 0$.  
Let $\sigma$ be a partition for $W$ and $\vec S$ be the associated skeleton graph. Note that $\vec S$ is symmetric. Let $S$ be the undirected graph obtained from $\vec S$  
by ignoring the orientations of the self-loops and 
by replacing every pair of oppositely oriented edges $\{u_iu_j, u_ju_i\}$, for $u_i\neq u_j$, with the undirected edge $(u_i, u_j)$. 

\begin{definition}
Let $f_1, \cdots, f_\ell$ be the edges of $S$. To each $f_j$, we associate the {\bf node-edge incidence vector}  $z'_j := \sum_{u_i\in f_j} e_i$. The {\bf node-edge incidence matrix} of $S$ is given by
$Z':= 
\begin{bmatrix}
z'_1 & \cdots & z'_\ell
\end{bmatrix}$.
\end{definition}

Let $\cX'$ be the convex cone spanned by $z'_1,\ldots, z'_\ell$. 
We establish the following result: 

\begin{lemma}\label{lem:xequalx'}
It holds that $\cX' = \cX$.
\end{lemma}

\begin{proof}
Note that each edge of $S$ corresponds to a cycle of $\vec S$; indeed, a self-loop $(u_i,u_i)$ corresponds to $u_iu_i$ and an edge $(u_i,u_j)$ between two distinct nodes corresponds to the $2$-cycle $u_iu_ju_i$.    
Relabel the cycles of $\vec S$ such that the first $\ell$ cycles $\vec C_j$, for $j = 1,\ldots, \ell$, correspond to the edges $f_j$ of $S$. It is clear that the node-cycle incidence vector $z_j$ of $\vec S$ coincides with the node-edge incidence vector $z'_j$ of $S$. Thus, $\cX' \subseteq \cX$. 
It remains to show that for any cycle $\vec C_j$ of $\vec S$, with length greater than $2$, the associated vector $z_j$ can be expressed as a nonnegative combination of the $z'_j$'s. We write $\vec C_j = u_1 u_2 \cdots u_\ell u_1$ for $\ell > 2$. Then, $f_1:= (u_1, u_2), f_2:= (u_2, u_3), \cdots, f_\ell := (u_\ell, u_1)$ are edges of $S$. It follows that 
$z_j = \frac{1}{2} \sum_{i = 1}^\ell z'_i$. 
\end{proof}

An immediate consequence of Lemma~\ref{lem:xequalx'} is that $\crk(Z) = \crk(Z')$. 
Also note that a symmetric digraph $\vec S$ is strongly connected if and only if $S$ is connected. 
The following result is thus a corollary of Theorem~\ref{thm:main} specializing to the class of graphons with symmetric support (which include the class of symmetric graphons). 

\begin{corollary}\label{cor:symmetricgraphon}
    Let $W$ be a step-graphon with symmetric support, and $\sigma$ be a partition for $W$. Let $x^*$ and $S$ be
    the associated concentration vector and the undirected skeleton graph. Further, let $Z'$ be the node-edge incidence matrix of $S$ and $\cX'$ be the convex cone spanned by the columns of $Z'$. Then, the following items hold:   
    \begin{enumerate}
    \item If $\crk(Z') > 0$ or if $x^*\notin \cX'$, then
    $$
    \lim_{n\to\infty} \bP(\vec G_n \sim W \mbox{ has a cycle cover}) = 0.
    $$
    \item If $\crk(Z') = 0$ and $x^*\in \rint \cX'$, and if $S$ is not connected, then 
    $$
     \lim_{n\to\infty} \bP(\vec G_n \sim W \mbox{ has a cycle cover}) = 1,
    $$
    and 
    $$
     \lim_{n\to\infty} \bP(\vec G_n \sim W \mbox{ has a Hamilton cycle}) = 0.
    $$
    \item If $\crk(Z') = 0$, $x^*\in \rint \cX'$, and $S$ is connected, then 
    $$
     \lim_{n\to\infty} \bP(\vec G_n \sim W \mbox{ has a Hamilton cycle}) = 1. 
    $$ 
    \end{enumerate}
\end{corollary}

In the earlier work~\cite{belabbas2021h,belabbas2023geometric}, we have addressed the weak $H$-property for the class of symmetric graphons $W$. However, the sampling procedure there is slightly different from the one used in this paper, as we elaborate below. 
In~\cite{belabbas2021h,belabbas2023geometric}, we first sample an undirected graph $G_n\sim W$, and then obtain the symmetric digraph $\vec G_n^\rs$ from $G_n$ by replacing each undirected edge with  a pair of oppositely oriented edges.  
To avoid any confusion, we say that $W$ has  {\em weak (resp., strong) $H^\rs$-property} if $\vec G_n^\rs\sim W$ has a cycle cover (resp., Hamilton cycle) \aas.  We establish the following result:

\begin{lemma}
A symmetric step-graphon $W$ has  weak (resp., strong) $H^\rs$-property if and only if it has weak (resp., strong) $H$-property.
\end{lemma}

\begin{proof}
Let $\overline W$ be the {\it saturation} of $W$, i.e., 
$$
\overline W(s,t) := 
\begin{cases}
1 & \mbox{if } W(s,t) \neq 0, \\
0 & \mbox{if } W(s,t) = 0.
\end{cases}
$$
It is clear that $\overline W$ and $W$ share the same support. By Corollary~\ref{cor:symmetricgraphon} (resp., the results of~\cite{belabbas2021h,belabbas2023geometric}), $\overline W$ has the weak/strong $H$-property (resp., $H^\rs$-property) if and only if $W$ does. 
It thus remains to show that $\overline W$ has the weak (resp., strong) $H$-property if and only if it has the weak (resp., strong) $H^\rs$-property. But this follows from the fact that the two sampling procedures,  $\vec G_n \sim \overline W$ and $\vec G^\rs_n \sim \overline W$, are equivalent with each other. To wit, since $\overline W$ takes value $1$ over its support, the two digraphs are completely determined by their respective empirical concentration vectors and, moreover, if $x(\vec G_n) = x(\vec G_n^\rs) =: x$, then $\vec G_n = \vec G_n^\rs = \vec K_{nx}$. We conclude the proof by pointing out that $x(\vec G_n)$ and $x(\vec G_n^\rs)$ are identically distributed.   
\end{proof}

\section{Proof of Proposition~\ref{prop:superregular}}\label{sec:blowuplemma}
The proof relies on the use of the Chernoff bound for Binomial random variable, which we recall below:

\begin{lemma}\label{lem:chernoffbound}
Suppose that $X\sim \mathrm{Bin}(N,p)$; then, for any $r \in [0,1]$, 
$$
\bP(X \leq (1 - r) Np ) \leq \exp\left ( - \frac{r^2}{2} Np\right ).
$$
\end{lemma}

Now, let $G_n\sim W^{\rs}$. Recall that $x(G_n)$ is the empirical concentration vector, which converges to $x^*$ \aas. 
It follows that \aas 
\begin{equation}\label{eq:regular1}
x_i(G_n)>  \frac{1}{2}\min \{x_i^* \mid i = 1,\ldots, m\} = : \alpha, \quad \mbox{for all } i = 1,\ldots, m.
\end{equation}
In the sequel, we assume that~\eqref{eq:regular1} holds. 
Let $\cE_{n}(u_i, u_j)$ be the event that the pair $(\pi^{-1}(u_i), \pi^{-1}(u_j))$ is $(\epsilon,\delta)$-super-regular. We establish the following result:

\begin{lemma}\label{lem:degreecondition}
For any $(u_i,u_j)\in E(S)$, the event $\cE_n(u_i,u_j)$ holds \aas.
\end{lemma}

\begin{proof}
For convenience, let $A:= \pi^{-1}(u_i)$ and $B:= \pi^{-1}(u_j)$. 
We show below that~\eqref{eq:edgedensitycondition} and~\eqref{eq:degreecondition} hold \aas.  

\xc{Proof that~\eqref{eq:edgedensitycondition} holds \aas.} 
For any given $X\subseteq A$ and $Y\subseteq B$, $e(X, Y)$ is a binomial $(|X||Y|, q_{ij})$ random variable. If 
\begin{equation}\label{eq:xyab}
|X| > \epsilon |A| \quad \mbox{and}  \quad |Y| > \epsilon |B|,  
\end{equation} 
then
\begin{multline*}
\bP(e(X,Y) \leq \delta |X||Y|)  = \bP\left (e(X, Y) \leq \left (1 - \frac{q_{ij} - \delta}{q_{ij}} \right ) q_{ij} |X||Y| \right )  \\
 \leq \exp\left ( - \frac{(q_{ij} - \delta)^2 |X||Y|}{2 q_{ij}}\right ) 
\leq \exp \left ( - \frac{(q_{ij} - \delta)^2 \epsilon^2 \alpha^2 }{2 q_{ij}} n^2 \right )  \leq \exp \left ( - \frac{q_{ij} \epsilon^2 \alpha^2 }{8} n^2 \right ),
\end{multline*}
where the first inequality follows from Lemma~\ref{lem:chernoffbound}, the second inequality follows from~\eqref{eq:regular1} and~\eqref{eq:xyab}, and the last inequality follows from the fact that $\delta \leq q_{ij}/2$ (see~\eqref{eq:defdelta}). 
The number of pairs that satisfy~\eqref{eq:xyab} is bounded above by the total number of pairs $(X, Y)\in 2^A \times 2^B$, which is $2^{|A| + |B|}\leq 2^n$. It follows that
$$
\bP(\mbox{event~\eqref{eq:edgedensitycondition} does {\em not} hold} ) \leq 2^n \exp \left ( - \frac{q_{ij} \epsilon^2 \alpha^2 }{8} n^2  \right ) \xrightarrow{n\to\infty} 0.
$$

\xc{Proof that~\eqref{eq:degreecondition} holds \aas.} For any $a\in A$, $e(a, B)$ is a  binomial $(|B|, q_{ij})$ random variable. Using the same arguments as above, we obtain that  
\begin{multline*}
\bP(e(a, B) \leq \delta |B|) = \bP\left ( e(a, B) \leq \left (1 - \frac{q_{ij} - \delta}{q_{ij}} \right ) q_{ij} |B| \right ) \\
\leq \exp\left ( - \frac{(q_{ij} - \delta)^2 |B|}{2 q_{ij}}\right ) 
\leq \exp \left (  - \frac{q_{ij} \alpha }{8} n \right ). 
\end{multline*}
Similarly, 
$$
\bP(e(b, A) \leq \delta |A|) \leq \exp \left (  - \frac{q_{ij} \alpha }{8} n \right ).
$$
We conclude that  
$$
\bP(\mbox{event~\eqref{eq:degreecondition} does {\em not} hold} )  \leq (|A| + |B|) \exp \left (  - \frac{q_{ij} \alpha }{8} n \right ) 
 \leq n \exp \left (  - \frac{q_{ij} \alpha }{8} n \right ) \xrightarrow{n\to\infty} 0. 
$$
This completes the proof of the lemma.
\end{proof}

Proposition~\ref{prop:superregular} is then an immediate consequence of Lemma~\ref{lem:degreecondition}; indeed, 
$$
\bP(G_n \sim W^\rs \mbox{ is } \epsilon\mbox{-}\delta\mbox{-super-regular}) \geq 1 - \sum_{(u_i,u_j)\in E(S)} \bP(\neg\cE_n(u_i,u_j)) \xrightarrow{n\to\infty} 1.
$$
This completes the proof. \hfill{\qed}
\end{document}